\documentclass{amsart}
\usepackage{amssymb, a4wide, mathdots, url,hyperref,graphicx}
\usepackage[usenames]{color}
\usepackage{enumerate}
\usepackage{wasysym}
\usepackage{dsfont}
\usepackage{tikz-cd}
\usepackage{comment}
\usepackage{MnSymbol}
\usepackage{eqnarray, amsmath}
\usepackage{kotex}
\def\N{{\mathbb N}} \def\Z{{\mathbb Z}}

 \def\Q{{\mathbb Q}}

\def\R{{\mathbb R}} \def\C{{\mathbb C}}
\def\A{{\mathbb A}}

\def\F{{\mathbb F}}

\def\lcm{\operatorname{lcm}}

\def\n{{\bf n}}

\newcommand{\bE}{\mathbb E}

\newcommand{\Nm}{\operatorname{Nm}}
\newcommand{\Spec}{\operatorname{Spec}}
\newcommand{\Div}{\operatorname{Div}}
\newcommand{\Sub}{\operatorname{Sub}}
\newcommand{\Supp}{{\rm Supp}}
\newcommand{\cl}{\operatorname{cl}}
\newcommand{\Gal}{\operatorname{Gal}}
\newcommand{\GL}{\operatorname{GL}}
\newcommand{\SL}{\operatorname{SL}}
\newcommand{\tr}{\operatorname{tr}}
\newcommand{\Frob}{\operatorname{Frob}}

\DeclareMathSymbol{\lhook}{\mathrel}{letters}{"2C}

\newtheorem{theorem}{Theorem}[section]
\newtheorem*{theorem*}{Theorem}
\newtheorem{lemma}[theorem]{Lemma}
\newtheorem{definition}[theorem]{Definition}
\newtheorem{proposition}[theorem]{Proposition}

\theoremstyle{remark}
\newtheorem{remark}[theorem]{Remark}
\newtheorem{hypothesis}{Hypothesis}

\title[Local vanishing density of order $1$ and Divisor moments]{Local vanishing density of order $1$ and Divisor moments}

\author{Wonwoong Lee}
\address{Department of
 Mathematics Education, Chonnam National University, Yongbong-ro, Buk-gu, Gwangju 61186, Chonnam Province, South Korea}
\email{dldnjsdnd041@gmail.com}

\thanks{}

\date{\today}



\begin{document}

\begin{abstract}
    We establish a general framework relating local vanishing densities of order $1$ to logarithmic growth of divisor moments on regular arithmetic curves. As applications, we recover classical divisor moment estimates for polynomial values and establish their function field analogue. Furthermore, under GRH, we obtain a result on divisor moments for polynomials in the Fourier coefficients of non-CM cuspidal eigenforms with trivial nebentypus and trivial inner twist. The Fourier coefficients need not be integers, and the polynomials need not be squarefree. Under additional analytic hypotheses, we also allow a Sato--Tate angle restriction.
\end{abstract}

\maketitle

\tableofcontents

\section{Introduction and Main results}\label{sec_Intro geometric}

\subsection{Main result}

The purpose of this paper is to rigorously establish a reasonable framework for the following naive principle, inspired by the joint work of Lau and the author \cite{LL26} and the author's previous work \cite{Lee26}:
\begin{align*}
    \text{Local vanishing density of order $1$} \quad \Longrightarrow \quad \text{Log power divisor moment}.
\end{align*}

Moreover, this framework will be developed from a geometric point of view, rather than within the classical setting adopted in the relevant papers \cite{Chi22}, \cite{LL26}, and \cite{Lee26}. This geometric perspective allows us to formulate our framework in a unified manner for both the arithmetic and function field cases. Since the framework includes a number field generalization of classical results, the assumption that the Fourier coefficients of modular forms are integers is no longer required in the results of \cite{GM14}, \cite{Chi22}, \cite{LL26}, and \cite{Lee26}. In particular, all the previous results will be recovered in the special case $X=\Spec \Z$ of our setting. We also provide the function field analogue of the divisor moment estimation for polynomial values. 

We emphasize an important point. In the previous works \cite{LL26} and \cite{Lee26}, the required mean value estimates for certain multiplicative functions associated with the $\lcm$-function were obtained using the multivariable Tauberian theorem. By contrast, the arguments established in this paper bypass this step under suitable $L^2$ and tail bound conditions. Accordingly, the results of this paper do not arise merely from recasting the earlier auxiliary results in geometric language. In particular, the applications of our main results require neither multivariable mean value estimates nor the multivariable Tauberian theorem.


\

 We introduce the following notations.

\begin{itemize}
    \item $X$: a normal integral separated $1$-dimensional scheme of finite type over $\Spec \Z$.
    \item $K$: the function field of $X$.
    \item $\Div(X)_{\geq 0}$: the set of effective divisors on $X$.
    \item $\Sub_X(D)$: the set of effective subdivisors of $D$ on $X$, i.e.,
    \begin{align*}
        \Sub_X(D):=\{D' \in \Div(X)_{\geq 0}: D-D' \geq 0\}.
    \end{align*}
    \item $\Supp(D)$: the support of $D$, i.e., if $D=\sum v_i Y_i$ and $Y_i$'s are prime divisors, then $$\Supp(D):=\bigcup_i Y_i.$$
    \item $\cl(X)$: the set of closed points of $X$, which is countable as $X$ is of finite type over $\Spec \Z$.
\end{itemize} 

A $1$-dimensional integral separated scheme of finite type over $\Spec \Z$ is called an {\it arithmetic curve}.\footnote{There are several conventions for this terminology. Some authors require regularity or flatness.} As $X$ is locally Noetherian and normal of dimension $\leq 1$, it is regular (see \cite[Lemma~28.12.6]{Stacks}). Thus $X$ is a regular arithmetic curve. The function field $K$ is a global field. If $K$ is a number field, then $X=\Spec O_{K,S}$, where $S$ is a finite set of places of $K$. If $K$ is a function field of characteristic $p>0$, then $X$ is an open subscheme of the unique connected smooth proper curve $C$ over $\F_q$ with function field $K$ (cf. \cite[p.~420]{Mor23}).

Let $q$ be a fixed rational prime. Let $\{A_x\}_{x\in (0,\infty)}$ be a collection of finite sets such that $A_i \subset A_j$ for $i\le j$ and $A:= \bigcup_x A_x$ is an infinite set. For a divisor-valued function $\n:A \to \Div(X)_{\geq 0}$ and an effective divisor $D \in \Div(X)_{\geq 0}$, let
    \begin{align*}
        \pi(x,D):=\#\{a \in A_x: \n(a)-D \geq 0\}.
    \end{align*}
Note that any prime divisor of $X$ can be seen as a closed point, since $X$ is $1$-dimensional and irreducible. Express an effective divisor $D$ as
\begin{align*}
    D=\sum_{\mathfrak l \in \cl(X)} v_{\mathfrak l}(D)[\mathfrak l].
\end{align*}

Let $\kappa(\mathfrak l)$ be the residue field at $\mathfrak l \in X$. For any closed point $\mathfrak l \in \cl(X)$, the residue field $\kappa(\mathfrak l)$ is finite. Indeed, since $X$ is of finite type over $\Spec \Z$, $\mathfrak l$ can be seen as a maximal ideal of a finite type $\Z$-algebra. It lies over a rational prime $(\ell) \subset \Z$ and the field extension $\kappa(\mathfrak l)/\kappa(\ell)$ is finite (cf. \cite[Proposition~10.35.19]{Stacks}). Thus, we can define the {\it absolute norm} of $D \in \Div(X)_{\geq 0}$ by
\begin{align*}
    \Nm D:=\prod_{\mathfrak l \in \cl(X)}|\kappa(\mathfrak l)|^{v_{\mathfrak l}(D)}.
\end{align*}
It is clear that the {\it logarithmic norm} $\log \Nm$ is completely additive. In other words, we have
\begin{align*}
    \Nm(D+D')=\Nm D\cdot \Nm D'
\end{align*}
for any $D,D' \in \Div(X)_{\geq 0}$. By definition, $\Nm[\mathfrak l]=|\kappa(\mathfrak l)|$ for a prime divisor $[\mathfrak l]$. The number of $\mathfrak l \in \cl(X)$ with $|\kappa(\mathfrak l)| \ll 1$ is finite (see Lemma~\ref{lem_prime number theorem geometric} below), so the number of effective divisors $D$ with $\Nm D \ll 1$ is also finite.

We say that $\n$ is polynomially bounded if there exists $\beta>0$ such that
\begin{align*}
    \Nm\n(a)<x^{\beta} \quad \text{if } a \in A_x
\end{align*}
for every $x \in (0,\infty)$.

\begin{definition}
    A function $h:\R_{>0}\times \Div(X)_{\geq 0} \to \R$ defined by
    \begin{align*}
        h(x,D):=\frac{\pi(x,D)}{|A_x|}
    \end{align*}
    is the local density function for the vanishing of $\n$.
\end{definition}

\begin{definition}\label{def_essentially multiplicative geometric}
    A function $g:\Div(X)_{\geq 0} \to \C$ is said to be essentially multiplicative if $g(0)=1$ and for any $D \in \Div(X)_{\geq 0}$, we have
    \begin{align*}
        g(D) \ll \prod_{v_{\mathfrak l}(D)>0} g(v_{\mathfrak l}(D)\mathfrak l)
    \end{align*}
    and further if every prime component $\mathfrak l$ of $D$ satisfies $|\kappa(\mathfrak l)| > B$ for some fixed constant $B>0$, we also have
    \begin{align*}
        g(D) \gg \prod_{v_{\mathfrak l}(D)>0} g(v_{\mathfrak l}(D)\mathfrak l).    \end{align*}
\end{definition}

We say that $h$ has the {\it sieve dimension sequence} $(\kappa_v)_{v\geq 1}$ if there exists $c \in (0,1)$ such that for a given integer $v \geq 1$,
\begin{align}\label{eqn_sieve dimension sequence estimation}
    \sum_{|\kappa(\mathfrak l)| \leq y}h(x,v[\mathfrak l]) =\kappa_v \log\log y +O(1)
\end{align}
uniformly in $x,y>0$ that satisfy $y^v \leq x^c$. If such a sequence $(\kappa_v)_{v\geq 1}$ exists, then we say that $h$ is of order $1$. Define
\begin{align}\label{eqn_Km definition}
    K_m:=\sum_{v \geq 1}\kappa_v((v+1)^m-v^m).
\end{align}

\begin{theorem}\label{thm_main_geometric}
    Let $X$ be a normal integral separated $1$-dimensional scheme of finite type over $\Spec \Z$. Let $\n: A \to \Div(X)_{\geq 0}$ be a polynomially bounded function and let $h$ be the local density function for the vanishing of $\n$. Suppose that there exists $c \in (0,1)$ for which $h(x,D)$ is essentially multiplicative in $D$ if $\Nm D \leq x^c$, and $h$ is of order $1$ with the sieve dimension sequence $(\kappa_v)_{v \geq 1}$. Assume that there is an integer $E \geq 1$ such that:
    \begin{enumerate}
        \item[\normalfont(i)] If $v>E$, then $\kappa_{v}=0$.
        \item[\normalfont(ii)] The function $h(x,v[\mathfrak l])$ satisfies
            \begin{align*}
            \sup_x\sum_{\substack{\mathfrak l \in \cl(X) \\ |\kappa(\mathfrak l)|^v \leq x^c}} h(x,v[\mathfrak l])^2 < \infty
            \end{align*} 
        if $1\leq v \leq E$.
        \item[\normalfont(iii)] There exists $\delta_0>0$ such that
        \begin{align*}
            \sup_x\sum_{\mathfrak l \in \cl(X)} \sum_{\substack{v>E \\ |\kappa(\mathfrak l)|^v \leq x^c}}|\kappa(\mathfrak l)|^{\delta_0 v}h(x,v[\mathfrak l]) < \infty.
        \end{align*}
    \end{enumerate}
    Then, for a positive integer $m$, we have
    \begin{align*}
        \sum_{a\in  A_x}
        |\Sub_X(\n(a))|^m \asymp |A_x|(\log x)^{K_m}.
    \end{align*}
\end{theorem}

Note that due to the polynomial bound condition, (ii) is equivalent to
\begin{align*}
    \sup_x\sum_{\substack{\mathfrak l \in \cl(X)}} h(x,v[\mathfrak l])^2 < \infty.
\end{align*}
Indeed, for a fixed $v$, we have
\begin{align*}
    \sum_{\substack{\mathfrak l \in \cl(X) \\ |\kappa(\mathfrak l)|^v > x^c}} h(x,v[\mathfrak l])^2 \leq \sum_{\substack{\mathfrak l \in \cl(X) \\ |\kappa(\mathfrak l)|^v > x^c}} h(x,v[\mathfrak l]) = \frac{1}{|A_x|}\sum_{a \in A_x} \#\{\mathfrak l \in \cl(X): v[\mathfrak l] \leq \n(a), \ |\kappa(\mathfrak l)|^v > x^c\}.
\end{align*}
Let $r=\#\{\mathfrak l \in \cl(X): v[\mathfrak l] \leq \n(a), \ |\kappa(\mathfrak l)|^v > x^c\}$. Then
\begin{align*}
    x^{cr} < \prod_{\substack{v[\mathfrak l] \leq \n(a) \\ |\kappa(\mathfrak l)|^v > x^c}} |\kappa(\mathfrak l)|^v \leq \Nm \n(a) < x^{\beta},
\end{align*}
so we have
\begin{align*}
    \frac{1}{|A_x|}\sum_{a \in A_x} \#\{\mathfrak l \in \cl(X): v[\mathfrak l] \leq \n(a), \ |\kappa(\mathfrak l)|^v > x^c\} < \frac{\beta}{c}.
\end{align*}

One sufficient condition for Theorem~\ref{thm_main_geometric}(iii) is that there exists a constant $\alpha>1/(E+1)$ such that
\begin{align}\label{eqn_h(x,vp) for v>E bound}
    h(x,v[\mathfrak l]) \ll v^B |\kappa(\mathfrak l)|^{-\alpha v}, \qquad \text{for all $v>E$ with $|\kappa(\mathfrak l)|^v \leq x^c$}.
\end{align}
In this case, we can simply take $\delta_0$ to be a positive real number that is smaller than $\alpha-1/(E+1)$.



\subsection{Applications}





We denote the average of values $g(a)$ over $a \in M$ by
\begin{align*}
    \bE_M[g(\cdot)]=\frac{1}{|M|}\sum_{a \in M}g(a).
\end{align*}
The first application of Theorem~\ref{thm_main_geometric} is the following classical result.

\begin{theorem}\label{thm_application 1: polynomial value}
    Let $F \in \Z[T]$ be a polynomial with integer coefficients. Denote by $\N(x;F)$ the set of positive integers $n$ with $n\leq x$ and $F(n) \neq 0$. Suppose that
    \begin{align*}
        F=cF_1^{e_1}F_2^{e_2}\cdot \cdots \cdot F_s^{e_s}
    \end{align*}
    for $c \in \Z$, primitive irreducible polynomials $F_1,\ldots,F_s$ that are pairwise distinct and $e_1, \ldots, e_s \geq 1$. Let $\kappa_v$ be the number of indices $1 \leq i\leq s$ for which $e_i \geq v$. Then,
    \begin{align*}
        \bE_{\N(x;F)}[d(F(\cdot))^m] \asymp_{F,m} (\log x)^{K_m},
    \end{align*}
    where $K_m$ is defined as in \eqref{eqn_Km definition}.
\end{theorem}

Note that
\begin{align*}
    K_m=\sum_{i=1}^s ((e_i+1)^m-1).
\end{align*}
When $F$ is irreducible, we have $\kappa_1=1$ and $\kappa_v=0$ for $v\geq 2$, so $K_m=2^m-1$. Theorem~\ref{thm_application 1: polynomial value} recovers Erd\"os's and Delmer's results \cite{Erd52}, \cite{Del71}.

The second application is its modular analogue, which generalizes \cite[Theorem~1.1]{Lee26} in two directions: the polynomial $F$ is not necessarily squarefree, and the Fourier coefficients of the modular form $f$ are not necessarily integers. For a number field $K$, let $\mathcal I_K$ be the set of all integral ideals of $\mathcal O_K$, the ring of integers of $K$. Define the divisor function $d_K: \mathcal I_K \to \N$ by letting $d_K(\mathfrak o)$ be the number of integral ideals dividing $\mathfrak o$. More explicitly,
\begin{align*}
    d_K(\mathfrak o)=(v_1+1)\cdot \cdots \cdot (v_r+1),
\end{align*}
where $v_1,\ldots,v_r$ are given by the prime ideal decomposition
\begin{align*}
    \mathfrak o=\mathfrak l_1^{v_1}\cdot \cdots \cdot \mathfrak l_r^{v_r}.
\end{align*}

We recall some basic facts on modular forms. By Deligne's bound, we have a bound $|a_f(p)| \leq 2p^{(k-1)/2}$ for every prime $p$. For this reason, we may denote by $\theta_f(p) \in [0,\pi]$ the Sato--Tate angle of $f$ at $p$, which is defined by
    \begin{align*}
        \frac{a_f(p)}{p^{(k-1)/2}}=2\cos \theta_f(p).
    \end{align*}

An inner twist of $f$ is a pair $(\gamma,\chi)$ of a $\Q$-embedding $\gamma: K_f \hookrightarrow \C$ and a Dirichlet character $\chi$ such that
\begin{align*}
    f^{\gamma} = f\otimes \chi.
\end{align*}
If $f$ is non-CM, then $\gamma \in {\rm Aut}(K_f/\Q)$ and $\chi$ is determined by $\gamma$ so that we may write $\chi=\chi_{\gamma}$. The set of all inner twists of $f$ forms a group $\Gamma$ in the natural way. We set
\begin{align*}
    H:=\bigcap_{\gamma \in \Gamma} \ker \chi_{\gamma},
\end{align*}
which is an open subgroup of $\Gal_{\Q}$. If $\Gamma=\{1\}$ (equivalently, $K^{\Gamma}=K$), then $H=\Gal_{\Q}$ obviously. Conversely, if $H=\Gal_{\Q}$, then we have $\chi_{\gamma}=1$ for all $\gamma \in \Gamma$, so
\begin{align*}
    f^{\gamma}=f.
\end{align*}
This means that $\gamma$ fixes every Fourier coefficient of $f$, so $\gamma=1$. Therefore, $H=\Gal_{\Q}$ is equivalent to $\Gamma=\{1\}$. In this case, we say that $f$ has trivial inner twist.

If $f$ is non-CM and has trivial nebentypus, then for any $\gamma \in \Gamma$, $\chi$ is either trivial or quadratic. Indeed, the nebentypus of $f \otimes \chi$ is $\chi^2$ and the nebentypus of $f^{\gamma}$ is $1^{\gamma}=1$, so we have $\chi^2=1$.
Thus one may consider the following sufficient condition which implies that a non-CM cuspidal eigenform $f$ of weight $k \geq 2$ and level $N$ with trivial nebentypus has trivial inner twist: either
\begin{enumerate}
    \item[\normalfont(a)] $[K_f:\Q]$ is odd, or
    \item[\normalfont(b)] $N$ is $1$ or squarefree.
\end{enumerate}
Suppose (a) holds. For any $\gamma \in \Gamma$ and $p \gg_f 1$, we have $$\gamma^2(a_f(p))=\gamma(\chi_{\gamma}(p)a_f(p))=\chi_{\gamma}(p)\gamma(a_f(p))=\chi_{\gamma}(p)^2a_f(p)=a_f(p). $$
Since $[K_f:\Q]$ is odd and so is $|{\rm Aut}(K_f/\Q)|$, this implies that $\gamma=1$. If $(b)$ holds, then having trivial inner twist follows from \cite[Theorem~3.9]{Rib80}. We remark that, although $f$ is assumed to be of weight~$2$ in \cite{Rib80}, the proof works for arbitrary weight $k\geq 2$.

For an integral ideal $\mathfrak o \subset \mathcal O_{K_f}$, let 
\begin{align*}
    \bar\rho_{f,\mathfrak o}: \Gal_{\Q} \longrightarrow \GL_2(\mathcal O_{K_f}/\mathfrak o)
\end{align*}
be the mod $\mathfrak o$ Galois representation associated with $f$.

\begin{hypothesis}\label{hypo_nice analytic properties of L-function}
    Let $K/\mathbb Q$ be a finite Galois extension with Galois group $G$. Suppose $G$ is isomorphic to the image of $\bar \rho_{f,\mathfrak o}$. Let $\chi$ be the character of a finite-dimensional irreducible complex representation $\rho_{\chi}$ of $\mathrm{Gal}(K/\mathbb Q)$. Then the $L$-function $L(s,\mathrm{Sym}^m f \otimes \chi)=L(s,\mathrm{Sym}^m \rho_f \otimes \rho_{\chi})$ satisfies the following properties.
    \begin{enumerate}
        \item[\normalfont(a)] The gamma factor of $L(s,\mathrm{Sym}^m f \otimes \chi)$ is
        \begin{align*}
            & L_{\infty}(s,\mathrm{Sym}^m f \otimes \chi)\\
            &=\begin{cases}
                \prod_{i=1}^{(m+1)/2}\Gamma_{\mathbb C}\left(s+\left(i-\frac{1}{2}\right)(k-1)\right)^{\chi(1)}, & \text{if $m$ is odd,} \\ \Gamma_{\mathbb R}(s+\epsilon)^{n^+(\chi)}\Gamma_{\mathbb R}(s+\epsilon')^{n^-(\chi)}\prod_{i=1}^{m/2}\Gamma_{\mathbb C}(s+i(k-1))^{\chi(1)}, &\text{if $m$ is even,}
            \end{cases}
        \end{align*}
        where $\Gamma_{\mathbb R}(s):=\pi^{-s/2}\Gamma(s/2)$, $\Gamma_{\mathbb C}(s):=\Gamma_{\mathbb R}(s)\Gamma_{\mathbb R}(s+1)$, and $\epsilon,\epsilon' \in \{0,1\}$, $\epsilon+\epsilon' \equiv 1 \pmod 2$, $\epsilon \equiv m/2 \pmod 2$.
        \item[\normalfont(b)] Define $\Lambda(s,\mathrm{Sym}^m f \otimes \chi):=L_{\infty}(s,\mathrm{Sym}^m f \otimes \chi)L(s,\mathrm{Sym}^m f \otimes \chi)$ and let $\delta(m,\chi)=1$ if $m=0$ and $\chi$ is trivial, and $\delta(m,\chi)=0$ otherwise. Then the function $$s^{\delta(m,\chi)}(1-s)^{\delta(m,\chi)}\Lambda(s,\mathrm{Sym}^m f \otimes \chi)$$ extends to an entire function on $\mathbb C$ of order $1$ which does not vanish at $s=0,1$.
        \item[\normalfont(c)] The function $\Lambda(s,\mathrm{Sym}^m f \otimes \chi)$ satisfies the functional equation
        \begin{align*}
             & \Lambda(s,\mathrm{Sym}^m f \otimes \chi)=\epsilon(\mathrm{Sym}^m f \otimes \chi)q(\mathrm{Sym}^m f\otimes \chi)^{1/2-s}\Lambda(1-s,\mathrm{Sym}^m f \otimes \overline{\chi}),
        \end{align*}
        where $\epsilon(\mathrm{Sym}^m f \otimes \chi)$ is a complex number of absolute value $1$ (known as the root number) and $q(\mathrm{Sym}^m f \otimes \chi)$ is the (arithmetic) conductor of $\mathrm{Sym}^m f \otimes \chi$.
        \item[\normalfont(d)] The analytic conductor $\mathfrak q(\mathrm{Sym}^m f \otimes \chi)$ satisfies the inequality
        \begin{align*}
            \mathfrak q(\mathrm{Sym}^m f\otimes \chi) \leq \mathfrak q(\mathrm{Sym}^m f)^{\chi(1)}\mathfrak q(\chi)^{m+1}. 
        \end{align*}
    \end{enumerate}
\end{hypothesis}

\begin{theorem}\label{thm_application 2: modular forms}
    Let $f$ be a non-CM cuspidal eigenform of weight $k\geq 2$ and level $N$ with trivial nebentypus. Denote by $K_f$ the number field associated with the Fourier coefficients of $f$. Suppose that $f$ has trivial inner twist (for instance, $[K_f:\Q]$ is odd, or $N$ is $1$ or squarefree). For a polynomial $F \in \mathcal O_{K_f}[T]$ and an interval $I \subseteq [0,\pi]$ of the positive Sato--Tate measure, let
    \begin{align*}
        \Pi(x,I;F,f):=\{p \leq x: p \nmid N, \ \theta_f(p) \in I, \ F(a_f(p)) \neq 0 \}.
    \end{align*}
    Assume GRH for the Dedekind zeta functions. If $I \neq [0,\pi]$, assume in addition Hypothesis~\ref{hypo_nice analytic properties of L-function} and GRH for such $L$-functions. Suppose that
    \begin{align*}
        F=cF_1^{e_1}F_2^{e_2}\cdot \cdots \cdot F_s^{e_s}
    \end{align*}
    for $c \in K_f^{\times}$, pairwise distinct primitive irreducible polynomials $F_1,\ldots,F_s \in K_f[T]$ and $e_1, \ldots, e_s \geq 1$. Let $\kappa_v$ be the number of indices $1 \leq i\leq s$ for which $v \leq e_i$. Then,
    \begin{align*}
        \bE_{\Pi(x,I;F,f)}[d_{K_f}(F(a_f(\cdot)))^m] \asymp_{F,f,m} (\log x)^{K_m},
    \end{align*}
    where $K_m$ is defined as in \eqref{eqn_Km definition}.
\end{theorem}

When $F$ is squarefree, we have $\kappa_1=s$ and $\kappa_v=0$ for $v \geq 2$, so $K_m=s(2^m-1)$. With $K_f=\Q$, this recovers \cite[Theorem~1.1]{Lee26}.

\begin{remark}
    Applying Theorem~\ref{thm_main_geometric} to get Theorem~\ref{thm_application 2: modular forms} is not as straightforward as one might expect, even under the assumption of trivial inner twists. There are two reasons for this. First, although our auxiliary divisor moment estimate (Theorem~\ref{thm_main_geometric}) works in a general framework, several additional lemmas need to be proved rigorously in order to use it to generalize \cite[Theorem~1.1]{Lee26} to the number field setting. In particular, such additional work is required when using the adelic open image theorem to show that the local vanishing density is essentially multiplicative. While the adelic open image theorem applies in the general number field setting, it yields only a decomposition of the local density function with respect to rational primes $\ell$. Further work is required to refine this to a decomposition with respect to the primes of $K_f$. These arguments are carried out in \S\S~\ref{subsec_non-CM cuspidal eigenform}.
Second, there is a gap in the proof of \cite[Theorem~A.3]{LL26}. More precisely, the part establishing the lower bound contains a gap in the treatment of the bad primes. Fortunately, this gap does not affect the overall arguments of either \cite{LL26} or \cite{Lee26}. However, filling this gap over a general number field is precisely what is needed for our purposes in this paper. Therefore, we give a correct treatment of the local vanishing densities at the bad primes. This is stated in Proposition~\ref{prop_local density mod forms bad prime essential multiplicativity}.
\end{remark}

We also consider applications to the function field analogue of these divisor moment results. Let $R=\F_q[U]$, $K=\F_q(U)$ and $F(T)=cF_1^{e_1}F_2^{e_2}\cdot \cdots \cdot F_s^{e_s} \in R[T]$, where $c \in R$ and $F_i$ are primitive polynomials in $R[T]$ which are irreducible over $K$. For a monic element $b$ in $R$ and for an effective divisor $D_b$ of the principal ideal $(b) \subseteq R$, we have $\Nm b=\Nm D_b=q^{\deg b}$. Note that
\begin{align*}
    |\Sub_X (F(D_b))|=\#\{d \in R: d \ \textrm{monic}, \ d \mid F(b)\} =: d_{\F_q[U]}(F(b))
\end{align*}
is the analogue of the number-of-divisors function. We define $\kappa_v$ as in Theorem~\ref{thm_application 1: polynomial value}, namely, as the number of indices $i$ with $e_i \geq v$.

\begin{theorem}\label{thm_application 3: function field polynomial_weak}
     Let
\begin{align*}
    S_{\leq x}:=\#\{b \in R: b \ \textrm{monic}, \ \Nm b \leq x,\ F(b) \neq 0\}. 
\end{align*}
Then we have
\begin{align}\label{eqn_polynomial divisor moment function field analogue weak ver}
    \mathbb E_{S_{\leq x}}[d_{\F_q[U]}(F(\cdot))^m]=\frac{1}{|S_{\leq x}|} \sum_{b \in S_{\leq x}} d_{\F_q[U]}(F(b))^m \asymp_{F,m,q} (\log x)^{K_m},  
\end{align}
where $K_m$ is defined as in \eqref{eqn_Km definition}.
\end{theorem}

 This is the direct function field analogue of Theorem~\ref{thm_application 1: polynomial value}. In fact, we can obtain a more natural and stronger statement as follows.
\begin{theorem}\footnote{ChatGPT suggested strengthening Theorem~\ref{thm_application 3: function field polynomial_weak} to averages over polynomials of fixed degree and provided the idea to obtain the present theorem using Theorem~\ref{thm_main_geometric}.}\label{thm_application 3: function field polynomial}
    Let
    \begin{align*}
        S_n:=\{b \in R: b \ \textrm{monic}, \ \deg b=n, \ F(b) \neq 0\}.
    \end{align*}
    Then,
    \begin{align*}
        \mathbb E_{S_n}[d_{\F_q[U]}(F(\cdot))^m] \asymp_{F,m,q} n^{K_m}.
    \end{align*}
    Consequently, \eqref{eqn_polynomial divisor moment function field analogue weak ver} holds.
\end{theorem}

We note that Wei and Yang's result \cite{WY25} provides an asymptotic formula of $$\mathbb E_{S_n}[d_{\F_q[T]}(F(b))]$$ for a quadratic polynomial of the form $F(U)=U(U+b_0)$, where $b_0 \in R$ is of degree $<n$.

Theorem~\ref{thm_application 3: function field polynomial} implies Theorem~\ref{thm_application 3: function field polynomial_weak} as follows. Put $n=\lfloor \log_q x \rfloor$. From the relation
\begin{align*}
    \sum_{b \in S_{\leq x}}d_{\F_q[U]}(F(b))^m = \sum_{r=0}^n \sum_{b \in S_r} d_R(F(b))^m,
\end{align*}
the lower bound is simply given by 
\begin{align*}
    \sum_{b \in S_n} d_{\F_q[U]}(F(b))^m \asymp_{F,m,q} |S_n|n^{K_m} \sim q^n n^{K_m}.
\end{align*}
The upper bound is also given by
\begin{align*}
    \sum_{r=0}^n q^r r^{K_m} \ll_{q} q^n n^{K_m}. 
\end{align*}
As $|S_{\leq x}|=(q^{n+1}-1)/(q-1)+O_F(1)$, these bounds imply Theorem~\ref{thm_application 3: function field polynomial_weak}.

\begin{remark}
    Although Theorem~\ref{thm_application 3: function field polynomial_weak} follows immediately from Theorem~\ref{thm_application 3: function field polynomial}, one can prove Theorem~\ref{thm_application 3: function field polynomial_weak} independently using Theorem~\ref{thm_main_geometric} with $A_x=S_{\leq x}$. This approach is  parallel to the argument of the proof of Theorem~\ref{thm_application 1: polynomial value}. 
\end{remark}

\subsection{Strategy of the proof}

For $D_1,\ldots, D_m \in \Div(X)_{\geq 0}$, we define the {\it join} of divisors $D_1,\ldots,D_m$ by
\begin{align*}
    D_1 \vee \cdots \vee D_m =\max\{D_1,\ldots,D_m\}:=\sum_{\mathfrak l \in \cl(X)} \max\{v_{\mathfrak l}(D_1),\ldots,v_{\mathfrak l}(D_m)\}[\mathfrak l].
\end{align*}
Since $D_1 \vee \cdots \vee D_m \leq D_1+\cdots+D_m$, we have
\begin{align}\label{eqn_norm inequality for join and sum divisor version}
    \Nm (D_1 \vee \cdots \vee D_m) \leq \Nm D_1\cdot \cdots \cdot \Nm D_m.
\end{align}

Note that
\begin{align*}
    \sum_{a \in A_x} |\Sub_X(\n(a))|^m&=\sum_{a \in A_ x}\bigg(\sum_{\substack {D_1 \in \Div(X)_{\geq 0} \\ D_1 \leq \n(a)}}1\bigg)\bigg(\sum_{\substack {D_2 \in \Div(X)_{\geq 0} \\ D_2 \leq \n(a)}}1\bigg)\cdot \cdots \cdot \bigg(\sum_{\substack {D_m \in \Div(X)_{\geq 0} \\ D_m \leq \n(a)}}1\bigg) \\
    &=\sum_{a \in A_x}\bigg(\sum_{\substack {D_1,\ldots,D_m \in \Div(X)_{\geq 0} \\ D_1,\ldots, D_m \leq \n(a)}}1\bigg)  \\
    &=\sum_{D_1, \ldots, D_m \in \Div(X)_{\geq 0}} \pi(x,D_1 \vee \cdots \vee D_m).
\end{align*}
Thus, we have
\begin{align*}
    \frac{1}{|A_x|}\sum_{a \in A_x} |\Sub_X(\n(a))|^m =\sum_{d_1,\ldots,d_m} h(x,D_1\vee \cdots \vee D_m).
\end{align*}

It will turn out that the main quantity we need to estimate is
\begin{align*}
    \sum_{\Nm D_1,\ldots,\Nm D_m \leq Y} h(x,D_1\vee \cdots \vee D_m)
\end{align*}
for an appropriate $Y$ depending on $x$. Define 
\begin{align*}
    \Phi_m(D):=\#\{(D_1,\ldots,D_m) \in \Div(X)_{\geq 0}^m: D_1\vee \cdots \vee D_m=D\}.
\end{align*}
By~\eqref{eqn_norm inequality for join and sum divisor version}, one has
\begin{align}\label{eqn_h(x,lcm) fundamental bound geometric version}
    \sum_{\Nm D \leq Y}\Phi_m(D)h(x,D) \leq \sum_{\Nm D_1,\ldots,\Nm D_m\leq Y} h(x,D_1\vee \cdots \vee D_m)\leq \sum_{\Nm D \leq Y^m}\Phi_m(D)h(x,D).
\end{align}

We will prove that the lower and upper bounds have the same log power. We will also prove that a tail part of the upper bound has the same log power, which yields the desired result.

\subsection{Notational Remark}

We use both notations $\# A$ and $|A|$ to denote the number of elements of a set $A$.

\subsection*{Acknowledgement}

The author has been supported by the National Research Foundation of Korea
(NRF) grant funded by the Korea government (MSIP) [No. RS-2024-00341327 and No. RS-2024-00415601
(G-BRL)].

\subsection*{AI disclosure}

The author was inspired by the possibility that existing results concerning divisor moment estimates could be obtained without the multivariable Tauberian theorem during interactions with OpenAI’s ChatGPT (5.6 Sol, accessed July 2026). This observation served as the starting point for the present paper. 
Following these initial interactions, mathematical frameworks, arguments, and proofs presented in this paper were independently established by the author, except where specific assistance from or suggestions by AI are explicitly acknowledged in the corresponding footnotes. The author also employed ChatGPT (6.0 Astra, accessed Oct 2026) for language editing and proofreading. The author takes full responsibility for the content of the paper.

\section{Arithmetic of logarithmic additive functions}

This section is devoted to developing some basic arithmetic properties of several functions related to effective divisors on $X$. All of these come from classical results in multiplicative number theory and their function field analogues.

\begin{lemma}\label{lem_log Phi_m additivity geometric}
    The function $\Phi_m$ is logarithmic additive. Namely, if $D_1, D_2 \in \Div(X)_{\geq 0}$ have no common prime component, then 
    \begin{align*}
        \log \Phi_m(D_1+D_2)=\log \Phi_m(D_1)+\log \Phi_m(D_2).
    \end{align*}
\end{lemma}

\begin{proof}
    We may assume $D_1=v_1[\mathfrak l_1]$ and $D_2=v_2[\mathfrak l_2]$. Then $\Phi_m(v_1[\mathfrak l_1]+v_2[\mathfrak l_2])$ is the number of $D_1,\ldots,D_m \in \Div(X)_{\geq 0}$ such that
    \begin{align*}
        D_i= u_{1,i}[\mathfrak l_1]+u_{2,i}[\mathfrak l_2], \qquad i=1,2,\ldots,m 
    \end{align*}
    with the maxima of $u_{1,i}$ and $u_{2,i}$ over $i$ equal to $v_1$ and $v_2$, respectively. This number equals $$\Phi_m(v_1[\mathfrak l_1])\Phi_m(v_2[\mathfrak l_2])$$ since the choice of $u_{1,i}$ and $u_{2,i}$ is independent.
\end{proof}

\begin{lemma}\label{lem_|Sub_X| is logarithmic subadditive}
    The function $|\Sub_X|$ is logarithmic additive and logarithmic subadditive. In other words, for any effective divisors $D_1, D_2 \in \Div(X)_{\geq 0}$, the following hold.
    \begin{enumerate}
        \item[\normalfont(1)] If $D_1$ and $D_2$ have no common divisor, we have
        \begin{align*}
            \log |\Sub_X(D_1+D_2)| = \log |\Sub_X(D_1)|+\log |\Sub_X(D_2)|
        \end{align*}
        \item[\normalfont(2)] We have
    \begin{align*}
        \log |\Sub_X(D_1+D_2)| \leq \log |\Sub_X(D_1)|+\log |\Sub_X(D_2)|.
    \end{align*}
    \end{enumerate}
\end{lemma}

\begin{proof}
    The first assertion follows immediately from the definition of $\Sub_X$. To prove the second assertion, write
    \begin{align*}
        D_1&=v_1[\mathfrak l_1]+\cdots +v_r[\mathfrak l_r], \\
        D_2&=v_1'[\mathfrak l_1]+\cdots +v_r'[\mathfrak l_r],
    \end{align*}
    with $v_1, v_1',\ldots, v_r,v_r' \geq 0$. Then we have
    \begin{align*}
        |\Sub_X(D_1+D_2)|&=(v_1+v_1'+1)\cdot \cdots \cdot  (v_r+v_r'+1) \\
        & \leq (v_1+1)(v_1'+1) \cdot \cdots \cdot (v_r+1)(v_r'+1).
    \end{align*}
    The last term equals $|\Sub_X(D_1)||\Sub_X(D_2)|$, which completes the proof.
\end{proof}

The following lemma is an analogue of the prime number theorem.

\begin{lemma}\label{lem_prime number theorem geometric}
    Let 
    \begin{align*}
        \pi_X(T):=\#\{\mathfrak l \in \cl(X): |\kappa(\mathfrak l)| \leq T\}.
    \end{align*}
    Then
    \begin{align*}
        \pi_X(T) \asymp_X \frac{T}{\log T}
    \end{align*}
    as $T \to \infty$.
\end{lemma}

\begin{proof}
    Suppose that ${\rm char}K(X)=0$. As we have mentioned in \S~\ref{sec_Intro geometric}, $X=\Spec O_{K,S}$ for some finite set of places $S$ of $K$. By adjoining finitely many points, we may assume that $X=\Spec O_K$, so in this case $\Div(X)_{\geq 0}=\mathcal I_K$ can be seen as the set of integral ideals of $K$ and
    \begin{align*}
        \pi_X(T)=\text{the number of prime ideals $\mathfrak l$ of $K$ with $\Nm \mathfrak l \leq T$.}
    \end{align*}
    The assertion follows from the Landau prime ideal theorem.

    Now suppose ${\rm char}K(X)=p>0$. As above, we may assume that $X$ is an affine curve with function field $K$. The set $\cl(X)$ can be seen be identified with the set of prime ideals of $\Gamma(X,O_X)$. According to \cite[Theorem~5.12]{Ros02}, we have
    \begin{align*}
        \#\{\mathfrak l \in \cl(X): \deg \mathfrak l=n\} \sim \frac{q^n}{n}.
    \end{align*}
    Note that $\deg \mathfrak l=[\kappa(\mathfrak l):\F_q]$ so that $|\kappa(\mathfrak l)|=q^{\deg \mathfrak l}$. We have
    \begin{align*}
        \pi_X(T)=\#\{\mathfrak l \in \cl(X): \deg \mathfrak l \leq \log_q T\} \sim q+\frac{q^2}{2}+\cdots+\frac{T}{\log_q T}. 
    \end{align*}
    Clearly, $\pi_X(T) \gg \frac{T}{\log T}$. The upper bound can also be easily obtained by
    \begin{align*}
        q+\frac{q^2}{2}+\cdots+\frac{q^{N}}{N}=\int_0^q \frac{t^N-1}{t-1}dt \ll \frac{q^N}{N}
    \end{align*}
    so that $q+\frac{q^2}{2}+\cdots+\frac{T}{\log_q T} \ll \frac{T}{\log T}$.
\end{proof}

\begin{lemma}\label{lem_omega bound geometric}
    Define a function $\omega_{\max,X}: \R_{>0} \to \N_0$ by
    \begin{align*}
        &\omega_{\max,X}(Y):=\max\{|\Supp(D)|:D \in \Div(X)_{\geq 0}, \ \Nm D \leq Y\}.
    \end{align*}
    Then
    \begin{align*}
        \omega_{\max,X}(Y) \asymp_X \frac{\log Y}{\log \log Y}.
    \end{align*}
\end{lemma}

\begin{proof}
    Let $r=|\Supp(D)|$ for $D\in \Div(X)_{\geq 0}$ with $\Nm D \leq Y$. Let $T=c_X r \log r$. According to Lemma~\ref{lem_prime number theorem geometric}, one may choose sufficiently small $c_X>0$ for which
    \begin{align*}
        \pi_X(T) \leq \frac{r}{2}.
    \end{align*}
    Thus, among the prime components of $D$, at least $r/2$ prime divisors have norm larger than $T$. Therefore,
    \begin{align*}
        \log \Nm D&=\sum_{\mathfrak l \in \cl(X)}v_{\mathfrak l}(D)\log |\kappa(\mathfrak l)| \geq \sum_{\substack{\mathfrak l \in \Supp(D) \\ |\kappa(\mathfrak l)|\geq T}}\log |\kappa(\mathfrak l)| \\
        & \geq \frac{r}{2}\log T \gg_X r\log r.
    \end{align*}
    This implies $r\log r \ll_X \log Y$ and hence $r \ll_X \log Y/\log\log Y$. This holds for every $r=|\Supp(D)|$ with $\Nm D \leq Y$, so $\omega_{\max,X}(Y) \ll_X \log Y/\log\log Y$.

    To obtain the lower bound, we define $D_T$ by
    \begin{align*}
        D_T:=\sum_{|\kappa(\mathfrak l)|\leq T}[\mathfrak l].
    \end{align*}
    By Lemma~\ref{lem_prime number theorem geometric}, the logarithmic norm of $D_T$ satisfies
    \begin{align*}
        \log \Nm D_T = \sum_{|\kappa(\mathfrak l)| \leq T}\log {|\kappa(\mathfrak l)|}  \leq \pi_X(T) \log T \ll_X T.
    \end{align*}
    If we take $T=c_X' \log Y$ for sufficiently small $c_X'>0$, then $\Nm D_T \leq Y$. For this $T$, again by Lemma~\ref{lem_prime number theorem geometric}, we have
    \begin{align*}
        |\Supp(D_T)|=\pi_X(T) \gg_X \frac{T}{\log T}\gg_X \frac{\log Y}{\log\log Y}
    \end{align*}
    and hence $\omega_{\max}(Y)\gg_X \log Y/\log\log Y$.
\end{proof}

Since $X$ is of finite type over $\Spec \Z$, the set $\cl(X)$ is countable. Thus we can enumerate $\cl(X)=\{\mathfrak l_i: i \in \N\}$ such that $|\kappa(\mathfrak l_i)| \leq |\kappa(\mathfrak l_j)|$ if $i \leq j$.

\begin{lemma}\label{lem_k and kappa p_k bound}
    There is a constant $C_X>0$ such that
    \begin{align*}
        k \leq C_X\frac{|\kappa(\mathfrak l_k)|}{\log |\kappa(\mathfrak l_k)|}
    \end{align*}
    for any $k \geq 1$.
\end{lemma}

\begin{proof}
    This immediately follows from Lemma~\ref{lem_prime number theorem geometric}, since $k \leq \pi_X(|\kappa(\mathfrak l_k)|)$.
\end{proof}

\begin{lemma}\label{lem_divisor bound geometric}
    We have
    \begin{align*}
        \limsup_{\Nm D \to \infty} \frac{(\log |\Sub_X(D)|)(\log\log \Nm D)}{\log \Nm D} \leq  2\log 2.
    \end{align*}
\end{lemma}

\begin{proof}
    We may extend the argument in \cite[Proposition~1.15]{Ove14} as follows. Suppose $\Nm D \geq 3$. Express $D$ as
    \begin{align*}
        D=v_1[\mathfrak l_1]+\cdots +v_r[\mathfrak l_r]
    \end{align*}
    with $|\kappa(\mathfrak l_1)| \leq \ldots \leq |\kappa(\mathfrak l_r)|$, $v_1, \ldots, v_r \geq 0$.  
    
    Since
    \begin{align*}
        \log \Nm D=v_1 \log|\kappa(\mathfrak l_1)| + \cdots + v_r\log |\kappa(\mathfrak l_r)|,
    \end{align*}
    we have
    \begin{align*}
        v_i \log 2 \leq v_i \log |\kappa(\mathfrak l_i)| \leq \log \Nm D
    \end{align*}
    for some $1 \leq i \leq r$ so that
    \begin{align}\label{eqn_log v_k+1 ineqaulity 1}
        \log(v_k+1) \leq \log \bigg(\frac{\log \Nm D}{\log 2}+1 \bigg).
    \end{align}
    We also have
    \begin{align}\label{eqn_log v_k+1 inequality 2}
        \log(v_i+1) \leq v_i \log 2
    \end{align}
    as $v_i \geq 0$.

    Let $\eta \in (1/2,1)$ and $\theta=\theta(D):=(1-\eta)\log\log \Nm D$. Choose the largest index $k$ for which $\log |\kappa(\mathfrak l_k)| \leq \theta$. By the logarithmic additivity of $|\Sub_X|$ (Lemma~\ref{lem_|Sub_X| is logarithmic subadditive}) and the bounds \eqref{eqn_log v_k+1 ineqaulity 1} and \eqref{eqn_log v_k+1 inequality 2}, we have
    \begin{align*}
        \log|\Sub_X(D)|
        &= \sum_{i=1}^k \log|\Sub_X(v_i[\mathfrak l_i])|+\sum_{i=k+1}^r \log|\Sub_X(v_i[\mathfrak l_i])| \\
        &= \sum_{i=1}^k \log(v_i+1) + \sum_{i=k+1}^r \log(v_i+1) \\
        & \leq \sum_{i=1}^k \log\bigg(\frac{\log \Nm D}{\log 2}+1 \bigg) + \sum_{i=k+1}^r v_i \log 2.
    \end{align*}
    Lemma~\ref{lem_k and kappa p_k bound} implies that
    \begin{align*}
        k \leq C_X|\kappa(\mathfrak l_k)|^2 \leq C_X \exp(2\theta),
    \end{align*}
    and we have $\theta^{-1}\log|\kappa(\mathfrak l_i)|>1$ for $i >k$. Therefore,
    \begin{align*}
        &\sum_{i=1}^k \log\bigg(\frac{\log \Nm D}{\log 2}+1 \bigg) + \sum_{i=k+1}^r v_i \log 2 
        \leq C_x \exp(2\theta) \log\bigg(\frac{\log \Nm D}{\log 2}+1 \bigg) + \frac{\log 2}{\theta}\sum_{i=k+1}^r v_i \log |\kappa(\mathfrak l_i)| \\
        & \qquad \leq  C_x \exp(2\theta) \log\bigg(\frac{\log \Nm D}{\log 2}+1 \bigg) + \frac{(\log 2)(\log \Nm D)}{\theta} \\
        & \qquad = C_X(\log \Nm D)^{2(1-\eta)}\log\bigg(\frac{\log \Nm D}{\log 2}+1 \bigg)+(1-\eta)^{-1}\frac{(\log 2)(\log \Nm D)}{\log \log \Nm D}.
    \end{align*}
    For arbitrary $\epsilon>0$, by taking $\eta$ close to $1/2$, we have
    \begin{align*}
        \log |\Sub_X(D)| \leq (2+\epsilon) \frac{(\log 2)(\log \Nm D)}{\log \log \Nm D} \qquad \text{if } \Nm D \geq N(\epsilon)
    \end{align*}
    for some $N(\epsilon)>0$.
\end{proof}

\section{Main Lemmas}\label{sec_main lemmas geometric}

In this section, we prove the main lemmas in the proof of Theorem~\ref{thm_main_geometric}.

Let $Y\geq Z$ be positive numbers and $S(Y,Z)$ be the set of effective divisors $D \in \Div(X)_{\geq 0}$ with $\Nm D \leq Y$ that is $Z$-{\it smooth}, i.e., for every prime component $\mathfrak l$ of $D$ we have $|\kappa(\mathfrak l)| \leq Z$. When $X=\Spec \Z$, this notion is exactly coincide with the notion of classical $Z$-smooth number. We also denote by $S'(Y,(z,Z])$ the set of $Z$-smooth effective divisors $D$ with $\Nm D > Y$ such that for every prime component $\mathfrak l$ of $D$, we require $z<|\kappa(\mathfrak l)|$. 

By the essential multiplicativity of $h(x,\cdot)$, we have
    \begin{align*}
        h(x,D)\ll \prod_{\substack{v_{\mathfrak l}(D)>0}}h(x,v_{\mathfrak l}(D)\mathfrak l)
    \end{align*}
if $\Nm D \leq x^c$. We let
    \begin{align*}
        h^{\times}(x,D):=\prod_{\substack{v_{\mathfrak l}(D)>0}}h(x,v_{\mathfrak l}(D)\mathfrak l).
    \end{align*}
This function is multiplicative, i.e., if $D$ and $D'$ have no common prime component, then
\begin{align*}
    h^{\times}(x,D+D')=h^{\times}(x,D)h^{\times}(x,D').
\end{align*}

\begin{lemma}\label{lem_Phi h(x,delta) sum order geometric ver}\footnote{While proofreading, we discovered a gap in the proof. ChatGPT pointed out that this gap could be remedied by adapting the idea in the draft of this paper and suggested a corresponding argument. The present proof follows this suggestion.}
    Let $h:\mathbb R_{>0} \times \Div(X)_{\geq 0} \to \R_{\geq 0}$ be a function. Assume the conditions in Theorem~\ref{thm_main_geometric} for $h$. Then for $Y\leq x^{c/E}$, 
    \begin{align*}
        \sum_{\substack{D \in \Div(X)_{\geq 0} \\ \Nm D \leq Y}} \Phi_m(D)h(x,D) \asymp (\log Y)^{K_m}, \quad \text{for } Y \gg 1.
    \end{align*}
\end{lemma}

\begin{proof}

    For $y>B$, where $B$ is a positive constant defined as in Definition~\ref{def_essentially multiplicative geometric} for $h$, let
    \begin{align*}
        P_x(y):=\prod_{B<|\kappa(\mathfrak l)|\leq y} \left(1+\sum_{1\leq v\leq E}\Phi_m(v[\mathfrak l])h(x,v[\mathfrak l]) \right).
    \end{align*}
    Then
    \begin{align*}
        \log P_x(y)&=\sum_{B<|\kappa(\mathfrak l)|\leq y}\log \bigg(1+\sum_{1\leq v\leq E}\Phi_m(v[\mathfrak l])h(x,v[\mathfrak l]) \bigg) \\
        &=\sum_{B<|\kappa(\mathfrak l)|\leq y}\bigg( \sum_{1\leq v \leq E} \Phi_m(v[\mathfrak l])h(x,v[\mathfrak l]) + \mathcal E_x(\mathfrak l)\bigg).
    \end{align*}
    Note that
    \begin{align*}
        \sum_{B<|\kappa (\mathfrak l)| \leq y}\mathcal E_x(\mathfrak l)=O(1),
    \end{align*}
    since
    \begin{align*}
        \mathcal E_x(\mathfrak l) \ll \bigg(\sum_{1\leq v \leq E}\Phi_m(v[\mathfrak l])h(x,v[\mathfrak l])\bigg)^2 \leq \sum_{1\leq v\leq E} \Phi_m(v[\mathfrak l])^2\sum_{1\leq v\leq E} h(x,v[\mathfrak l])^2 \ll \sum_{1\leq v\leq E} h(x,v[\mathfrak l])^2,
    \end{align*}
    and hence
    \begin{align*}
        \sum_{\mathfrak l} \mathcal E_x(\mathfrak l) \ll \sup_x\sum_{1\leq v\leq E}\sum_{\substack{\mathfrak l \in \cl(X) \\ |\kappa(\mathfrak l)|^v\leq x^c}} h(x,v[\mathfrak l])^2=O(1)
    \end{align*}
    by condition (ii) of Theorem~\ref{thm_main_geometric}. Consequently,
    \begin{align*}
        \sum_{B< |\kappa (\mathfrak l| \leq y}\bigg(\sum_{1\leq v \leq E}\Phi_m(v[\mathfrak l])h(x,v[\mathfrak l])+\mathcal E_x(\mathfrak l)\bigg) &=\sum_{B<|\kappa (\mathfrak l)| \leq y}\sum_{1\leq v\leq E}((v+1)^m-v^m)h(x,v[\mathfrak l])+O(1) \\
        &=\sum_{1\leq v\leq E}((v+1)^m-v^m)\sum_{B< |\kappa (\mathfrak l)|\leq y}h(x,v[\mathfrak l])+O(1) \\
        &=\sum_{1\leq v\leq E}((v+1)^m-v^m)\kappa_v \log\log y+O(1) \\
        &=K_m \log\log y+O(1),
    \end{align*}
    so we obtain
    \begin{align*}
        P_x(y) \asymp (\log y)^{K_m}.
    \end{align*}

    Define
    \begin{align*}
        T_x'(Y,y):=\sum_{\substack{D\in S'(Y,(B,y]) \\  \forall \mathfrak l, \ 0\leq v_{\mathfrak l}(D) \leq E}} \Phi_m(D)h^{\times}(x,D).
    \end{align*}
    We use the extension of Rankin's trick in the geometric framework. For any $\theta>0$, we have
    \begin{align*}
        T_x'(Y,y) &\leq Y^{-\theta} \sum_{\substack{D\in S'(Y,(B,y]) \\  \forall \mathfrak l, \ 0\leq v_{\mathfrak l}(D) \leq E}} \Phi_m(D)h^{\times}(x,D) (\Nm D)^{\theta} \\
        & \leq Y^{-\theta}\prod_{B<|\kappa(\mathfrak l)|\leq y}\bigg(1+\sum_{1\leq v \leq E}((v+1)^m-v^m)h(x,v[\mathfrak l])|\kappa(\mathfrak l)|^{v\theta}\bigg).
     \end{align*}
     One has
     \begin{align*}
         \frac{T_x'(Y,y)}{P_x(y)} \leq Y^{-\theta} \prod_{B<|\kappa(\mathfrak l)| \leq y} \frac{1+\sum_{1\leq v \leq E}((v+1)^m-v^m)h(x,v[\mathfrak l])|\kappa(\mathfrak l)|^{v\theta}}{1+\sum_{1\leq v\leq E}\Phi_m(v[\mathfrak l])h(x,v[\mathfrak l]) }
     \end{align*}
     and consequently
     \begin{align*}
         \log T_x'(Y,y) - \log P_x(y) \leq -\theta\log Y + \sum_{B<|\kappa(\mathfrak l)| \leq y}  \log \bigg(\frac{1+\sum_{1\leq v \leq E}((v+1)^m-v^m)h(x,v[\mathfrak l])|\kappa(\mathfrak l)|^{v\theta}}{1+\sum_{1\leq v\leq E}\Phi_m(v[\mathfrak l])h(x,v[\mathfrak l]) }\bigg).
     \end{align*}
     Note that
     \begin{align*}
         \log \frac{1+x_1+x_2}{1+x_1} \leq x_2
     \end{align*}
     for any $x_1,x_2 \geq 0$. Indeed, this follows easily from
     \begin{align*}
         1+x_1+x_2 \leq e^{x_2}(1+x_1), \qquad \forall x_1,x_2 \geq 0.
     \end{align*}
     Thus we have
     \begin{align*}
         \log T_x'(Y,y) - \log P_x(y) &\leq -\theta\log Y + \sum_{B<|\kappa(\mathfrak l)| \leq y} \sum_{1 \leq v\leq E}((v+1)^m-v^m)h(x,v[\mathfrak l])(|\kappa(\mathfrak l)|^{v\theta}-1)
     \end{align*}
     Let $\alpha \in (0,1)$ and choose $y=Y^{\alpha}$ and $\theta=1/(E\log y)$. Let $\mathfrak h_{x,v}(i):=\sum_{|\kappa(\mathfrak l)| = i}h(x,v[\mathfrak l])$ and $H_{x,v}(t):=\sum_{|\kappa(\mathfrak l)| \leq t}h(x,v[\mathfrak l])=\sum_{i \leq t} \mathfrak h_{x,v}(i)$. Then 
    \begin{align*}
        \sum_{B<|\kappa(\mathfrak l)| \leq y}h(x,v[\mathfrak l])(|\kappa(\mathfrak l)|^{v\theta}-1)=\sum_{B<i\leq y} \mathfrak h_{x,v}(i) (i^{v\theta}-1).
    \end{align*}
    Since $y^E<Y^E \leq x^c$, we may apply \eqref{eqn_sieve dimension sequence estimation}. More precisely, by partial summation, if $1\leq v \leq E$ we have
    \begin{align*}
        \sum_{B<i\leq y} \mathfrak h_{x,v}(i) (i^{v\theta}-1)&\leq (y^{v\theta}-1)H_{x,v}(y)-v\theta\int_B^y H_{x,v}(t)t^{v\theta-1}dt \\
        &=\kappa_v (y^{v\theta}-1)\log\log y -v\theta\int_B^y \kappa_v t^{v\theta-1}\log\log t \ dt + O(y^{v\theta}) \\
        &=\kappa_v (y^{v\theta}-1)\log\log y-\kappa_v \bigg(t^{v\theta}\log\log t \bigg|_B^y -\int_B^y t^{v\theta}\frac{1}{t\log t} dt\bigg) + O(y^{v\theta})\\
        &=\kappa_v \int_B^y \frac{t^{v\theta-1}}{\log t}dt-\kappa_v\log\log y+ O(y^{v\theta}).
    \end{align*}
    Note that $y^{v\theta} \leq y^{E\theta}=y^{1/\log y}=e$. In addition, we have
    \begin{align*}
        \int_B^y \frac{t^{v\theta-1}}{\log t}dt=\int_{\log B}^{\log y} \frac{e^{xv\theta}}{x}dx 
        &=\log\log y-\log\log B+\int_{\log B}^{\log y}\frac{e^{xv\theta}-1}{x}dx \\
        &=\log\log y+\int_{v\theta\log B}^{v\theta\log y}\frac{e^w-1}{w}dw+O(1) \\
        &=\log\log y+O(y^{v\theta}), 
    \end{align*}
    so
    \begin{align*}
        \kappa_v \int_B^y \frac{t^{v\theta-1}}{\log t}dt-\kappa_v\log\log y+ O(y^{v\theta})=O(1).
    \end{align*}
    This implies that
    \begin{align*}
        \log T_x'(Y,y)-\log P_x(y) &\leq -\theta \log Y + O(1) \\
        & = -\frac{1}{\alpha E}+O(1).
    \end{align*}
    We may choose $\alpha$ sufficiently close to $0$ so that 
    \begin{align*}
        \log T_x'(Y,y)-\log P_x(y) < -\log 2,
    \end{align*}
    which is equivalent to
    \begin{align*}
        P_x(y)-T_x'(Y,y)>\frac{1}{2}P_x(y).
    \end{align*}
    The lower bound follows from
    \begin{align*}
        \sum_{\substack{D \in \Div(X)_{\geq 0} \\ \Nm D \leq Y}} \Phi_m(D)h(x,D) \gg P_x(y)-T_x'(Y,y) \gg P_x(y) \asymp (\log y)^{K_m}.
    \end{align*}

    We can get the upper bound similarly. We have
    \begin{align*}
         \log \bigg( \sum_{\substack{D \in \Div(X)_{\geq 0} \\ \Nm D \leq Y}} \Phi_m(D)h(x, D)\bigg) 
         &\leq \log \bigg(\prod_{|\kappa (\mathfrak l)| \leq Y}\sum_{\substack{v \geq 0 \\ |\kappa(\mathfrak l)|^v \leq Y}}\Phi_m(v[\mathfrak l])h(x,v[\mathfrak l])\bigg) +O(1) \\
         &\leq \sum_{|\kappa (\mathfrak l)| \leq Y}\bigg(\sum_{1\leq v \leq E}\Phi_m(v[\mathfrak l])h(x,v[\mathfrak l])+\mathcal E_x(\mathfrak l)\bigg)+\sum_{\mathfrak l \in \cl(X)} R_x(\mathfrak l)+O(1),
    \end{align*}
    where
    \begin{align*}
        R_x(\mathfrak l):=\sum_{\substack{v>E \\ |\kappa(\mathfrak l)|^v \leq Y}}\Phi_m(v[\mathfrak l])h(x,v[\mathfrak l]).
    \end{align*}
    Clearly, one has
    \begin{align*}
        \sum_{\mathfrak l \in \cl(X)} R_x(\mathfrak l) \ll \sum_{\mathfrak l \in \cl(X)} \sum_{\substack{v>E \\ |\kappa(\mathfrak l)|^v \leq x^c}} |\kappa (\mathfrak l)|^{\delta_0 v}h(x,v[\mathfrak l]) =O(1),
    \end{align*}
    so it is negligible. 
    
    For $1\leq v\leq E$, since $|\kappa (\mathfrak l)|^v\leq Y^E \leq x^c$, we can use the same argument as for the lower bound to get
    \begin{align*}
         \sum_{|\kappa (\mathfrak l)| \leq Y}\bigg(\sum_{1\leq v \leq E}\Phi_m(v[\mathfrak l])h(x,v[\mathfrak l])+\mathcal E_x(\mathfrak l)\bigg)=\log(\log Y)^{K_m}+O(1),
    \end{align*}
    which completes the proof.
\end{proof}

\begin{lemma}\label{lem_smooth tail bound lemma 1 geometric}
    Let $\eta \in (0,1)$ and $h:\R_{>0} \times \Div(X)_{\geq 0} \to \R_{\geq 0}$ be a function. Assume the conditions in Theorem~\ref{thm_main_geometric} for $h$. Let $0<Y\leq x^c$ and $Z>0$. For arbitrary $\theta \in (0,\delta_0)$ (which may depend on $Y$ and $Z$), we have
    \begin{align}\label{eqn_smooth tail bound geometric}
        \log\bigg(\sum_{\substack{D \in S(Y,Z) \\ \Nm D \geq Y^{\eta}}}\Phi_m(D)h(x,D)\bigg) \leq K_m\log\log Z+O(Z^{E\theta})-\frac{\eta\theta}{2}\log Y.
    \end{align}
\end{lemma}

\begin{proof}

     To show \eqref{eqn_smooth tail bound geometric}, we will bound the logarithm of
    \begin{align*}
        T_x(Y,Z)=\sum_{\substack{D \in S(Y,Z) \\ \Nm D \geq Y^{\eta}}}\Phi_m(D)h^{\times}(x,D)
    \end{align*}
    using the extension of Rankin's trick. One has
    \begin{align*}
        T_x(Y,Z) &\leq Y^{-\eta\theta/2}\sum_{\substack{D \in S(Y,Z) \\ \Nm D \geq Y^{\eta}}}\Phi_m(D)h^{\times}(x,D)(\Nm D)^{\theta/2} \\
        &\leq Y^{-\eta\theta/2}\prod_{\substack{|\kappa(\mathfrak l)| \leq Z}}\bigg(1+\sum_{\substack{v\geq 1 \\ |\kappa(\mathfrak l)|^v \leq Y}}((v+1)^m-v^m)h(x,v[\mathfrak l])|\kappa(\mathfrak l)|^{v\theta/2}\bigg),
    \end{align*}
    so
    \begin{align*}
        \log T_x(Y,Z) &\leq -\frac{\eta \theta}{2} \log Y+\sum_{|\kappa(\mathfrak l)| \leq Z}\log\bigg(1+\sum_{\substack{v\geq 1 \\ |\kappa(\mathfrak l)|^v \leq Y}}((v+1)^m-v^m)h(x,v[\mathfrak l])|\kappa(\mathfrak l)|^{v\theta/2}\bigg) \\
        &\leq -\frac{\eta \theta}{2} \log Y +\sum_{v\geq 1}((v+1)^m-v^m)\sum_{|\kappa(\mathfrak l)| \leq \min\{Z,Y^{1/v}\}}h(x,v[\mathfrak l])|\kappa(\mathfrak l)|^{v\theta/2} \\
        &=-\frac{\eta \theta}{2} \log Y + \sum_{1 \leq v\leq E}((v+1)^m-v^m)\sum_{|\kappa(\mathfrak l)| \leq \min\{Z,Y^{1/v}\}}h(x,v[\mathfrak l])|\kappa(\mathfrak l)|^{v\theta} \\
        & \quad+ \sum_{v>E}((v+1)^m-v^m)\sum_{|\kappa(\mathfrak l)| \leq \min\{Z,Y^{1/v}\}}h(x,v[\mathfrak l])|\kappa(\mathfrak l)|^{v\theta}.
    \end{align*}
    As in the proof of Lemma~\ref{lem_Phi h(x,delta) sum order geometric ver}, we have
    \begin{align*}
        \sum_{|\kappa(\mathfrak l)| \leq \min\{Z,Y^{1/v}\}}h(x,v[\mathfrak l])|\kappa(\mathfrak l)|^{v\theta} \leq \kappa_v \int_2^Z \frac{t^{v\theta-1}}{\log t}dt+ O(Z^{v\theta}).
    \end{align*}
    Thus,
    \begin{align*}
        \sum_{1 \leq v\leq E}\sum_{|\kappa(\mathfrak l)|\leq \min\{Z,Y^{1/v}\}}h(x,v[\mathfrak l])|\kappa(\mathfrak l)|^{v\theta} \leq \sum_{1\leq v \leq E}\left(\kappa_v \int_2^Z \frac{t^{v\theta-1}}{\log t}dt+O(Z^{v\theta})\right).
    \end{align*}
    As in the proof of Lemma~\ref{lem_Phi h(x,delta) sum order geometric ver}, we have
    \begin{align*}
        \int_2^Z \frac{t^{v\theta-1}}{\log t}dt=\log\log Z+O(Z^{v\theta}). 
    \end{align*}
    
    It remains to estimate the sum on $v> E$. In fact, this part is negligible, since
    \begin{align*}
        &\sum_{v>E}((v+1)^m-v^m) \sum_{|\kappa(\mathfrak l)| \leq \min\{Z,Y^{1/v}\}}h(x,v[\mathfrak l])|\kappa(\mathfrak l)|^{v\theta} \\ &\ll_{m,\delta_0} \sum_{v>E}\sum_{|\kappa(\mathfrak l)|\leq \min\{Z,Y^{1/v}\}}h(x,v[\mathfrak l])|\kappa(\mathfrak l)|^{v\delta_0} =O(1)
    \end{align*}
    according to condition (iii) in Theorem~\ref{thm_main_geometric}.
\end{proof}

We apply Lemma~\ref{lem_smooth tail bound lemma 1 geometric} with specific choices of $\theta$ and $Z$ depending on $Y$. Fix $\lambda \in (0,1/E)$ and $\theta_0\in (0,\delta_0)$. Let
\begin{align*}
    u:=\frac{\log Y}{\log Z}, \qquad \theta:=\min\bigg\{\lambda\frac{\log u}{\log Z}, \theta_0\bigg\}.
\end{align*}
From the definition, we have
\begin{align*}
    -\frac{\eta\theta}{2}\log Y=-\frac{\eta \lambda u}{2}\log u, \qquad Z^{E\theta}=u^{E\lambda},
\end{align*}
when $\lambda\frac{\log u}{\log Z} \leq \theta_0$ (equivalently, $u\log u\leq \theta_0\lambda^{-1} \log Y$). Choose $u_0>0$ depending on $\eta$ and $\lambda$ for which the big $O$-term in Lemma~\ref{lem_smooth tail bound lemma 1 geometric} satisfies
\begin{align*}
    O(Z^{E\theta})\leq Cu^{E\lambda} \leq \frac{\eta\lambda}{4}u\log u \qquad \text{if } u\geq u_0.
\end{align*}
Then we have
\begin{align*}
    Cu^{E\lambda}-\frac{\eta\theta}{2}\log Y \leq \frac{\eta\lambda}{4} u\log u-\frac{\eta\lambda}{2}u\log u \leq -\frac{\eta\lambda}{4}u\log u.
\end{align*}
Therefore, if $u_0\leq u$ and $u\log u\leq \theta_0\lambda^{-1} \log Y$,
\begin{align*}
    \sum_{\substack{\delta \in S(Y,Y^{1/u}) \\ \Nm \delta \geq Y^{\eta}}}\Phi_m(\delta)h(x,\delta) \ll \bigg(\frac{1}{u}\log Y\bigg)^{K_m}\exp\bigg(-\frac{\eta\lambda}{4}u\log u\bigg).
\end{align*}

If $\lambda\frac{\log u}{\log Z}>\theta_0$ (equivalently, $u\log u > \theta_0\lambda^{-1} \log Y$), then $u>Z^{\theta_0/\lambda}>Z^{E\theta_0}$ so that
\begin{align*}
    \log Y=u\log Z>Z^{E\theta_0}\log Z,
\end{align*}
and hence 
\begin{align*}
    O(Z^{E\theta}) \leq \frac{\eta\theta_0}{4}\log Y.
\end{align*}
In this case, Lemma~\ref{lem_smooth tail bound lemma 1 geometric} implies 
\begin{align*}
    \sum_{\substack{D \in S(Y,Y^{1/u}) \\ \Nm D \geq Y^{\eta}}}\Phi_m(D)h(x,D) \ll (\log Y)^{K_m}Y^{-\frac{\eta\theta_0}{4}}.
\end{align*}

In short, we have the following.

\begin{lemma}\label{lem_smooth tail bound lemma 2 geometric}
    Let $\eta \in (0,1)$. If $Y \leq x^c$, there exist positive constants $u_0$, $\theta_0$, and $\lambda<1/E$ such that
    \begin{enumerate}
        \item[\normalfont(1)] If $u_0\leq u$ and $u\log u \leq \theta_0 \lambda^{-1} \log Y$, then
        \begin{align*}
            \sum_{\substack{D \in S(Y,Y^{1/u}) \\ \Nm D \geq Y^{\eta}}}\Phi_m(D)h(x,D) \ll \bigg(\frac{1}{u}\log Y\bigg)^{K_m}\exp\bigg(-\frac{\eta\lambda}{4}u\log u\bigg).
        \end{align*}
        \item[\normalfont(2)] If $u\log u > \theta_0 \lambda^{-1} \log Y$, then 
        \begin{align*}
            \sum_{\substack{D \in S(Y,Y^{1/u}) \\ \Nm D \geq Y^{\eta}}}\Phi_m(D)h(x,D) \ll_{\epsilon} Y^{-\frac{\eta\theta_0}{4}+\epsilon}.
        \end{align*}
    \end{enumerate}
\end{lemma}

\section{Proof of Theorem~\ref{thm_main_geometric}}

To prove Theorem~\ref{thm_main_geometric}, we use the lemmas in \S~\ref{sec_main lemmas geometric} and follow a variant of the arguments in \cite{LL26} and \cite{Lee26}, which originate in \cite{Erd52}.

\subsection{Upper bound}

Note that the contribution from $\n(a)=0$ is clearly negligible; since $\#\Sub_X(0)=1$, we have
\begin{align*}
    \sum_{\substack{a \in A_x, \\ \n(a)=0}}|\Sub_X(\n(a))|^m \leq |A_x|.
\end{align*}

Express
\begin{align*}
    \n(a)=[\mathfrak l_1]+  [\mathfrak l_2] + \cdots + [\mathfrak l_J],  
\end{align*}
where $[\mathfrak l_i]$ are prime components of $\n(a)$, not necessarily pairwise distinct, and $|\kappa(\mathfrak l_1)| \leq |\kappa(\mathfrak l_2)| \leq \ldots \leq |\kappa(\mathfrak l_J)|$. We let $P_i:=[\mathfrak l_1]+ [\mathfrak l_2] + \cdots + [\mathfrak l_i]$ and $j=j(\n(a))$ be the largest index for which $\Nm P_j\leq x^{c'}$, where $c':=c/mE$.

The sum of $|\Sub_X(\n(a))|^m$ over $a\in A_x$ with $j=0$ is also negligible; if $j=0$, then every prime component $[\mathfrak l]$ of $\n(a)$ satisfies $|\kappa(\mathfrak l)|>x^{c'}$, which implies $x^{Jc'}<\Nm \n(a)$. The polynomial bound condition for $\n$ implies $J<\beta/c'$, so
\begin{align*}
    |\Sub_X(\n(a))|^m\leq 2^{\beta m/c'}=O(1).
\end{align*}
Hereafter, we always assume $j\geq 1$.

Let $H:=(2\beta+1)/c'$. We divide the sum into $\Sigma_<+\Sigma_{\geq}$, where $\Sigma_<$ (resp. $\Sigma_{\geq}$) is the sum over $j$ with $J-j < H$ (resp. $J-j \geq H$). By Lemma~\ref{lem_|Sub_X| is logarithmic subadditive}, we have
\begin{align*}
    \sum_{\substack{a \in A_x \\ J-j<H}} |\Sub_X(\n(a))|^m 
    &\ll \sum_{\substack{a \in A_x \\ J-j<H}}|\Sub_X(P_j)|^m \\ 
    &= \sum_{\substack{a \in A_x \\ J-j<H}}\bigg(\sum_{\substack {D_1,\ldots,D_m \in \Div(X)_{\geq 0} \\ D_1,\ldots, D_m \leq P_j}}1\bigg) \\
    &\leq \sum_{\Nm D_1,\ldots,\Nm D_m \leq x^{c'}} \pi(x,D_1 \vee \cdots \vee D_m) \\
    &=\sum_{\Nm D_1,\ldots,\Nm D_m \leq x^{c'}} |A_x|h(x,D_1 \vee \cdots \vee D_m). 
\end{align*}
By~\eqref{eqn_h(x,lcm) fundamental bound geometric version} and Lemma~\ref{lem_Phi h(x,delta) sum order geometric ver}, it is bounded by
\begin{align*}
    \leq |A_x| \sum_{\Nm D \leq x^{c/E}}\Phi_m(D)h(x,D) \ll |A_x|(\log x)^{K_m}
\end{align*}
as $x \gg 1$. 

So we only need to bound $\Sigma_{\geq}$. For this sum, we may assume $J-j \geq H$. If $|\kappa(\mathfrak l_{j+1})| \geq x^{c'/2}$, then one has
\begin{align*}
    \Nm \n(a) \geq |\kappa(\mathfrak l_{j+1})||\kappa(\mathfrak l_{j+2})| \cdot \cdots \cdot |\kappa(\mathfrak l_{J})| \geq x^{c'H/2}=x^{\beta+\frac{1}{2}},
\end{align*}
which contradicts the assumption $\Nm \n(a) \leq x^{\beta}$. Thus we have $|\kappa(\mathfrak l_{j+1})| < x^{c'/2}$ and consequently, $\Nm P_j > x^{c'/2}$. Let $s=s(a)$ be the positive integer satisfying
\begin{align*}
    x^{\frac{1}{s+1}}\leq  |\kappa(\mathfrak l_j)| <x^{\frac{1}{s}}.
\end{align*}
Let $X:=x^{c'm}$ so that $x^{1/s}=X^{1/c'ms}$, $x^{c'/4}=X^{1/4m}$. We divide the sum into two cases: $c'ms \log(c'ms) \leq \theta_0 \lambda^{-1} \log X$ or not, where $\lambda$ and $\theta_0$ are chosen to satisfy Lemma~\ref{lem_smooth tail bound lemma 2 geometric} for $u=c'ms$. We denote the corresponding sums by $\Sigma_{\geq,1}$ and $\Sigma_{\geq,2}$, respectively. For the former case, note that since $|\Sub_X|$ is logarithmic subadditive (Lemma~\ref{lem_|Sub_X| is logarithmic subadditive}), we have
\begin{align*}
    |\Sub_X(\n(a))|^m & \leq |\Sub_X(P_j)|^m |\Sub_X(\mathfrak l_{j+1})|^m \cdot \cdots \cdot |\Sub_X(\mathfrak l_{J})|^m \\
    & \leq |\Sub_X(P_j)|^m 2^{m(J-j)} \\
    & \leq |\Sub_X(P_j)|^m 2^{m\beta(s+1)}.
\end{align*}
The last inequality follows from
\begin{align*}
    x^{\beta} \geq \Nm n(a) \geq |\kappa(\mathfrak l_{j+1})||\kappa(\mathfrak l_{j+2})| \cdot \cdots \cdot |\kappa(\mathfrak l_{J})| \geq x^{\frac{J-j}{s+1}}.
\end{align*}
Moreover, since we have $\log \Nm (D-D') + \log \Nm D' = \log \Nm D$ for $D' \in \Sub_X(D)$, at least one of $D-D'$ and $D$ has logarithmic norm at least $\frac{1}{2}\log \Nm D$. This implies that
\begin{align*}
    |\Sub_X(D)| \leq 2\#\{D' \in \Sub_X(D): \Nm D' \geq \sqrt{\Nm D}\}.
\end{align*}
Consequently, together with the fact that $\Nm P_j > x^{c'/2}$, we conclude that every effective subdivisor $D'$ of $P_j$ with $\Nm D' \geq \sqrt{\Nm P_j}$ satisfies $\Nm D' \geq x^{c'/4}$. Therefore, we obtain
\begin{align*}
    |\Sub_X(\n(a))|^m \ll 2^{m\beta s} |\Sub_X(P_j)|^m \ll 2^{m\beta s}\sum_{\substack{D_1,\ldots,D_m \in S(x^{c'},x^{1/s}) \\ D_1 \vee \cdots \vee D_m \leq \n(a) \\ \Nm D_1,\ldots,\Nm D_m \geq x^{c'/4}}} 1.
\end{align*}
Thus,
\begin{align*}
    \Sigma_{\geq,1} &\ll \sum_{c'ms \log(c'ms) \leq \theta_0 \lambda^{-1} \log X} 2^{m\beta s}\sum_{\substack{D_1,\ldots,D_m \in S(x^{c'},x^{1/s})\\ \Nm D_1,\ldots,\Nm D_m \geq x^{c'/4}}} \pi(x,D_1 \vee \cdots \vee D_m) \\
    &= |A_x|\sum_{c'ms \log(c'ms) \leq \theta_0 \lambda^{-1} \log X} 2^{m\beta s}\sum_{\substack{D_1,\ldots,D_m \in S(x^{c'},x^{1/s})\\ \Nm D_1,\ldots,\Nm D_m \geq x^{c'/4}}} h(x,D_1 \vee \cdots \vee D_m) \\
    &\ll |A_x|\sum_{c'ms \log(c'ms) \leq \theta_0 \lambda^{-1} \log X} 2^{m\beta s}\sum_{\substack{D \in S(x^{c'm},x^{1/s})\\ \Nm D \geq x^{c'/4}}}\Phi_m(D)h(x,D).
\end{align*}
The sum over $s$ with $c'ms \leq u_0$ for a constant $u_0$ is negligible; it is easily bounded by 
\begin{align*}
    \ll (\log x)^{K_m}
\end{align*}
according to Lemma~\ref{lem_Phi h(x,delta) sum order geometric ver}. Now we apply Lemma~\ref{lem_smooth tail bound lemma 2 geometric} for
\begin{align*}
    Y=X=x^{c'm}, \qquad X=x^{1/s}, \qquad u=c'ms, \qquad \eta=\frac{1}{4m}
\end{align*}
to obtain 
\begin{align*}
    &|A_x|\sum_{\substack{c'ms> u_0 \\ c'ms\log(c'ms) \leq \theta_0\lambda^{-1} \log X}} 2^{m\beta s}\sum_{\substack{D \in S(x^{c'm},x^{1/s})\\ \Nm D \geq x^{c'/4}}}\Phi_m(D)h(x,D)
    \\
    &\ll_{\epsilon}|A_x|\sum_{\substack{c'ms> u_0 \\ c'ms\log(c'ms) \leq \theta_0\lambda^{-1} \log X}} 2^{m\beta s}(\log X)^{K_m}\exp\bigg(-\frac{\eta\lambda}{4}c'ms\log (c'ms)\bigg).
\end{align*}
Since
\begin{align*}
    \sum_{\substack{c'ms>u_0 \\ c'ms\log(c'ms) \leq \theta_0\lambda^{-1} \log X}} 2^{m\beta s}\exp\bigg(-\frac{\eta\lambda}{4}c'ms\log (c'ms)\bigg) 
    &\ll 
    \sum_{s} 2^{m\beta s}\exp\bigg(-\frac{\eta\lambda}{4}c'ms\log (c'ms)\bigg) \\
    &\leq \sum_s \exp(Cs)\exp\bigg(-\frac{\eta\lambda}{4}c'ms\log (c'ms)\bigg) \\
    &=O(1)
\end{align*}
for some $C>0$, we have
\begin{align*}
    |A_x|\sum_{\substack{c'ms> u_0 \\ c'ms\log(c'ms) \leq \theta_0\lambda^{-1} \log X}} 2^{m\beta s}(\log X)^{K_m}\exp\bigg(-\frac{\eta\lambda}{4}c'ms\log (c'ms)\bigg) \ll |A_x|(\log x)^{K_m}.
\end{align*}

The remaining part is the sum over $a$ for which $c'ms\log(c'ms)>\theta_0\lambda^{-1} \log X$. Under such a condition, by combining the logarithmic subadditivity of $|\Sub_X|$, Lemma~\ref{lem_divisor bound geometric}, and the polynomial bound $|\kappa(\mathfrak l_i)| < x^{\beta}$, we get
\begin{align*}
    |\Sub_X(\n(a))|^m &\leq |\Sub_X(\mathfrak l_{j+1}+ \cdots + \mathfrak l_J)|^m |\Sub_X(P_j)|^m \ll \exp\bigg(\frac{C_{m,\beta}\beta \log x}{\log\log x}\bigg)\sum_{\substack{D_1,\ldots,D_m \in S(x^{c'},x^{1/s}) \\ D_1 \vee\cdots \vee D_m \leq \n(a) \\ \Nm D_1,\ldots,\Nm D_m \geq x^{c'/4}}} 1
\end{align*}
for some $C_{m,\beta}>0$. Since $2\leq |\kappa(\mathfrak l_j)| \leq x^{1/s}$ (by the definition of $s$), the condition $s \leq \frac{\log x}{\log 2}$ is already implied, and hence
\begin{align*}
    \Sigma_{\geq,2}
    &\ll \sum_{\substack{c'ms \log(c'ms) > \theta_0 \lambda^{-1} \log X \\ s \leq \frac{\log x}{\log 2}}} \sum_{\substack{a \in A_x \\ s(a)=s}} \exp\bigg(\frac{C_{m,\beta} \log x}{\log\log x}\bigg)\sum_{\substack{D_1,\ldots,D_m \in S(x^{c'},x^{1/s}) \\ D_1 \vee\cdots \vee D_m \leq \n(a) \\ \Nm D_1,\ldots,\Nm D_m \geq x^{c'/4}}} 1
    \\ 
    &\leq \sum_{\substack{c'ms \log(c'ms) > \theta_0 \lambda^{-1} \log X \\ s \leq \frac{\log x}{\log 2}}} \exp\bigg(\frac{C_{m,\beta} \log x}{\log\log x}\bigg)\sum_{\substack{D_1,\ldots,D_m \in S(x^{c'},x^{1/s})\\ \Nm D_1, \ldots, \Nm D_m \geq x^{c'/4}}} \pi(x,D_1 \vee \cdots \vee D_m) \\
    &= |A_x|\sum_{\substack{c'ms \log(c'ms) > \theta_0 \lambda^{-1} \log X \\ s \leq \frac{\log x}{\log 2}}} \exp\bigg(\frac{C_{m,\beta} \log x}{\log\log x}\bigg)\sum_{\substack{D_1,\ldots,D_m \in S(x^{c'},x^{1/s})\\ \Nm D_1, \ldots, \Nm D_m \geq x^{c'/4}}} h(x,D_1 \vee \cdots \vee D_m) \\
    &\leq |A_x|\sum_{\substack{c'ms \log(c'ms) > \theta_0 \lambda^{-1} \log X \\ s \leq \frac{\log x}{\log 2}}} \exp\bigg(\frac{C_{m,\beta} \log x}{\log\log x}\bigg)\sum_{\substack{D \in S(x^{c'm},x^{1/s})\\ \Nm D \geq x^{c'/4}}}\Phi_m(D)h(x,D)
\end{align*}
By Lemma~\ref{lem_smooth tail bound lemma 2 geometric}, this is bounded by
\begin{align*}
    \ll_{\epsilon} 
    |A_x|\sum_{\substack{c'ms \log(c'ms) > \theta_0 \lambda^{-1} \log X \\ s \leq \frac{\log x}{\log 2}}} \exp\bigg(\frac{C_{m,\beta} \log x}{\log\log x}\bigg) x^{-\frac{\eta\theta_0 c'm}{4}+\epsilon} = O(|A_x|).
\end{align*}

\subsection{Lower bound}

The proof of the lower bound is much simpler. We have
\begin{align*}
    \sum_{a \in A_{x}}|\Sub_X(n(a))|^m 
    &\geq \sum_{a \in A_{x}}\sum_{\substack{D_1,\ldots,D_m \leq n(a) \\ \Nm D_1, \ldots, \Nm D_m \leq x^{c/E}}}1 \\
    &= |A_x|\sum_{\Nm D_1,\ldots,\Nm D_m \leq x^{c/E}} h(x,D_1 \vee \cdots \vee D_m). 
\end{align*}
By~\eqref{eqn_h(x,lcm) fundamental bound geometric version} and Lemma~\ref{lem_Phi h(x,delta) sum order geometric ver}, the right-hand side of the above is bounded below by
\begin{align*}
    \geq |A_x|\sum_{\Nm D \leq x^{c/E}}\Phi_m(D)h(x,D) \asymp |A_x|(\log x)^{K_m}.
\end{align*}

\section{Applications to divisor moments}

\subsection{Polynomial values}\label{subsec_polynomial_values}

Let $X=\Spec \Z$. In this case, we can identify $\Div(X)_{\geq 0}$ with the set of ideals of $\Z$, or equivalently, the set $\N$. For $D \in \Div(X)_{\geq 0}$ (=$(\delta) \subset \Z$), the set of effective subdivisors of $D$ is identified with
\begin{align*}
    \Sub_X(D)=\{m \in \N: m \mid \delta\},
\end{align*}
so $|\Sub_X(D)|=d(\delta)$. Each closed point $\mathfrak l \in \cl(X)$ is a rational prime $\ell$, and its residue field is $\F_\ell$. We have $|\kappa(\mathfrak l)|=\ell$, and hence $\Nm D=\delta \in \N$.

\begin{proof}[Proof of Theorem~\ref{thm_application 1: polynomial value}]
    We set $A_x=\N(x;F)$ and $\n(n)=F(n) \in \N$ for $n \in A_x$ (more precisely, $\n(n)$ is the ideal generated by $F(n)$). The local density function $h(x,\delta)$ for the vanishing of $\n$ is given by
    \begin{align}\label{eqn_h(x,delta) and N_F(delta) polynomial}
        h(x,\delta)=\frac{\#\{n \leq x: \delta \mid F(n)\}}{x}=\frac{N_F(\delta)}{\delta}\bigg(1+O\bigg(\frac{\delta}{x}\bigg)\bigg),
    \end{align}
    where 
    \begin{align*}
        N_F(\delta):=\#\{a \bmod \delta: F(a) \equiv 0 \pmod \delta\}.
    \end{align*}
    Since $N_F(\delta)$ is multiplicative, $\delta \leq x^c$, and $\omega(\delta) \ll \log x$, we have
    \begin{align*}
        h(x,\delta)&=(1+o(1))\prod_{\ell^v \mid\mid \delta}h(x,\ell^v)
    \end{align*}
    so that $h$ is essentially multiplicative.
    
    For each $i$ and $\ell^v$, if we let 
    \begin{align*}
        h_i(x,\ell^v)=\frac{\#\{n \leq x: \ell^v \mid F_i(n)^{e_i}\}}{x},
    \end{align*}
    then for all but finitely many primes $\ell$, we have
    \begin{align*}
        h(x,\ell^v)=\sum_{i=1}^s h_i(x,\ell^v).
    \end{align*}
    It is obvious that
    \begin{align*}
        \ell^v \mid F_i(n)^{e_i} \ \text{if and only if} \ \ell^{\lceil v/e_i \rceil} \mid F_i(n). 
    \end{align*}
    Suppose $1 \leq v \leq e_i$. In this case, as $\lceil v/e_i\rceil=1$, we have
    \begin{align*}
        h_i(x,\ell^v)=\frac{\#\{n \leq x: \ell \mid F_i(n)\}}{x} = \frac{N_{F_i}(\ell)}{\ell}+O\bigg(\frac{1}{x}\bigg)
    \end{align*}
    so that
    \begin{align*}
        \sum_{\ell\leq y} h_i(x,\ell^v) = \sum_{\ell\leq y} \frac{N_{F_i}(\ell)}{\ell}+O\bigg(\frac{\pi(y)}{x}\bigg)=\kappa_v \log\log y + O(1)+ O\bigg(\frac{\pi(y)}{x}\bigg)
    \end{align*}
    (cf. \cite[(10)]{Kow12}). For $y^v \leq x^c$, we have $\pi(y)/x \ll 1.$ On the other hand, if $v>e_i$, or equivalently, $\lceil v/e_i \rceil \geq 2$, then we have
    \begin{align*}
        \#\{n \leq x: \ell^{\lceil v/e_i \rceil} \mid F_i(n)\} = \frac{N_{F_i}(\ell)x}{\ell^{\lceil v/e_i\rceil}}+O(1)
    \end{align*}
    and hence
    \begin{align*}
        \sum_{\ell \leq x}h_i(x,\ell^v) \leq \sum_{\ell \leq x} \bigg(\frac{\deg F_i}{\ell^2}+O\bigg(\frac{1}{x}\bigg)\bigg)=O(1).
    \end{align*}
    
    Therefore, for a given $v \geq 1$, we have
    \begin{align*}
        \sum_{\ell \leq x}h(x,\ell^v)=\kappa_v \log\log x +O(1)
    \end{align*}
    where $\kappa_v$ is the number of indices $i$ with $v \leq e_i$. In other words, $h$ has sieve dimension sequence $(\kappa_v)_{v \geq 1}$. Condition (iii) of Theorem~\ref{thm_main_geometric} follows from \eqref{eqn_h(x,vp) for v>E bound}, by taking $\delta_0^{-1}=2E(E+1)$. Condition~(ii) of Theorem~\ref{thm_main_geometric} also follows immediately, since $h(x,\ell^v) \ll \ell^{-1}$ for $1 \leq v \leq E$. We complete the proof by applying Theorem~\ref{thm_main_geometric}.
\end{proof}

\subsection{Non-CM cuspidal eigenforms}\label{subsec_non-CM cuspidal eigenform}

Next, we suppose that $X=\Spec \mathcal O_E$ for a number field $E$. In this case, $D \in \Div(X)_{\geq 0}$ is identified with an integral ideal $\mathfrak o \in \mathcal I_E$. Each closed point $\mathfrak l$ can be interpreted as a prime ideal of $E$, and the norm of $D$ equals the norm of $\mathfrak o$
\begin{align*}
    \Nm D=N_{E/\Q} (\mathfrak o).
\end{align*}
The number of effective subdivisors of $D$ is the same as $d_E(\mathfrak o)$. We assume that $f$ is a non-CM cuspidal eigenform satisfying the assumptions in Theorem~\ref{thm_application 2: modular forms}.

Let $E=K_f$. We set 
    \begin{align*}
        A_x=\Pi(x,I;F,f)_{\rm good}=\{p \in \Pi(x,I;F,f): p \gg_{F,f} 1\}
    \end{align*}
    and $\n(p)=F(a_f(p))\mathcal O_K \in \mathcal I_K$ for $p \in \Pi(x,I;F,f)$. One may easily verify the polynomial bound $\Nm \n(p) \ll x^{\beta}$ for $p \in A_x$ using Deligne's bound.
    The associated local density function $h$ is given by
    \begin{align*}
        h(x,\mathfrak o)=\frac{\#\{p \in \Pi(x,I;F,f)_{\rm good}: \mathfrak o \mid F(a_f(p))\}}{|\Pi(x,I;F,f)_{\rm good}|}.
    \end{align*}
    
    For a prime $\mathfrak l$ of $E$, let 
    \begin{align*}
        \rho_{f,\mathfrak l}: \Gal_{\Q} \longrightarrow \GL(V_{\mathfrak l})
    \end{align*}
    be a $2$-dimensional $\mathfrak l$-adic representation attached to $f$. Choosing a stable $\mathcal O_{E,\mathfrak l}$-lattice $T_{\mathfrak l} \subset V_{\mathfrak l}$ and an appropriate basis for $V_{\mathfrak l}$, we can view $\rho_{f,\mathfrak l}$ as
    \begin{align*}
        \rho_{f,\mathfrak l}: \Gal_{\Q} \longrightarrow \GL(T_{\mathfrak l}) \simeq \GL_2(\mathcal O_{E,\mathfrak l}).
    \end{align*}
    Its mod $\mathfrak l$ reduction gives
    \begin{align*}
        \bar\rho_{f,\mathfrak l}: \Gal_{\Q} \longrightarrow \GL_2(\kappa(\mathfrak l)),
    \end{align*}
    where $\kappa(\mathfrak l)$ is the residue field at $\mathfrak l$, and more generally, we have the mod $\mathfrak l^v$ reduction
    \begin{align*}
        \bar\rho_{f,\mathfrak l^v}: \Gal_{\Q} \longrightarrow \GL_2(\mathcal O_{E,\mathfrak l}/\mathfrak l^v \mathcal O_{E,\mathfrak l}).
    \end{align*}
    We also let
    \begin{align*}
        \rho_f: \Gal_{\Q} \longrightarrow \GL_2(\A_{E,fin}) \simeq \GL_2(E \otimes \A_{\Q,fin})
    \end{align*}
    be the adelic Galois representation attached to $\rho_f$ and $\bar \rho_{f,\mathfrak o}$ be its mod $\mathfrak o$ reduction
    \begin{align*}
        \bar \rho_{f,\mathfrak o}=\prod_{\mathfrak l^v  \mid\mid \mathfrak o}\bar \rho_{f,\mathfrak l^v}: \Gal_{\Q} \longrightarrow \GL_2(\mathcal O_E/ \mathfrak o) \simeq \prod_{\mathfrak l^v \mid\mid \mathfrak o} \GL_2(\mathcal O_{E,\mathfrak l^v}/\mathfrak l^v \mathcal O_{E,\mathfrak l^v}).
    \end{align*}
    Define
    \begin{align*}
        &\bar{\mathcal A}_{f,\mathfrak o}:= \bar \rho_{f,\mathfrak o}(\Gal_{\Q}), \\
        &\mathcal C_{F,f,\mathfrak o}:=\{\sigma \in \bar{\mathcal A}_{f,\mathfrak o}: F(\tr(\sigma))\equiv 0 \pmod{\mathfrak o}\},
    \end{align*}
    and for $d \in (\mathcal O_{E,\mathfrak l}/\mathfrak l^v \mathcal O_{E,\mathfrak l})^{\times}$, $t \in \mathcal O_{E,\mathfrak l}/\mathfrak l^v \mathcal O_{E,\mathfrak l}$,
    \begin{align*}
        &C_{\mathfrak l^v}(d,t):=\{g \in \GL_2(\mathcal O_{E,\mathfrak l}/\mathfrak l^v \mathcal O_{E,\mathfrak l}): \det g=d, \ \tr g = t\}, \\
        &A_{\mathfrak l^v}(d):=\{g \in \GL_2(\mathcal O_{E,\mathfrak l}/\mathfrak l^v \mathcal O_{E,\mathfrak l}): \det g=d\}.
    \end{align*}
    
We first extend \cite[Lemma~5.2]{Lee26} to the number field setting. Let $q_{\mathfrak l}:=|\kappa(\mathfrak l)|$ and
    \begin{align*}
        \frac{t^2}{4}-d=\pi_{\mathfrak l}^v u, \qquad u \in (\mathcal O_{E,\mathfrak l}/\mathfrak l ^{n-v}\mathcal O_{E,\mathfrak l})^{\times},
    \end{align*}
where $\pi_{\mathfrak l}$ is a uniformizer of $\mathfrak l$.

\begin{lemma}\label{lem_A lambda^n d, t number}
    The number of elements of $C_{\mathfrak l^n}(d,t)$ is $O(q_{\mathfrak l}^{2n})$. Moreover, if $\ell>2$, it is given by 
    \begin{align*}
        |C_{\mathfrak l^n}(d,t)|=\begin{cases}
            \sum_{j=0}^{n-1}(q_{\mathfrak l}-1)q_{\mathfrak l}^{2n-\lceil j/2 \rceil -1}+q_{\mathfrak l}^{n+\lfloor n/2 \rfloor} & \text{if } v=n, \\
            \sum_{j=0}^{v}(q_{\mathfrak l}-1)q_{\mathfrak l}^{2n-\lceil j/2 \rceil -1}+2q_{\mathfrak l}^{2n-v/2-1} & \text{if $v<n$, $v$ is even, and $u$ is a square mod $\mathfrak l$}, \\
            \sum_{j=0}^{v}(q_{\mathfrak l}-1)q_{\mathfrak l}^{2n-\lceil j/2 \rceil -1} & \text{otherwise},
        \end{cases}
    \end{align*}
    The number of elements of $A_{\mathfrak l^n}(d)$ is given by
    \begin{align*}
        |A_{\mathfrak l^n}(d)|=q_{\mathfrak l}^{3(n-1)}|\SL_2(\F_{q_{\mathfrak l}})|=q_{\mathfrak l}^{3n-2}(q_{\mathfrak l}^2-1).
    \end{align*}
    Consequently,
    \begin{align*}
        \frac{|C_{\mathfrak l^n}(d,t)|}{|A_{\mathfrak l^n}(d)|} = \frac{1}{q^n}\bigg(1+O\bigg(\frac{1}{q}\bigg)\bigg),
    \end{align*}
    where the big $O$-term is uniform in $n,d,t$.
\end{lemma}

\begin{proof}
    First, suppose $\ell>2$. The following argument for $C_{\mathfrak l^n}(d,t)$ is is an immediate generalization of the proof of \cite[Lemma~5.2]{Lee26}. Write $g= \begin{pmatrix}x & y \\ z & t-x \end{pmatrix} \in \GL_2(\mathcal O_{E,\mathfrak l}/\mathfrak l^n \mathcal O_{E,\mathfrak l})$. We have to compute the number of triples $(x,y,z) \in (\mathcal O_{E,\mathfrak l}/\mathfrak l^n \mathcal O_{E,\mathfrak l})^3$ for which $x(t-x)-yz \equiv d \pmod{\mathfrak l^n}$. This equation is equivalent to
    \begin{align*}
        \bigg(x-\frac{t}{2}\bigg)^2 \equiv \pi_{\mathfrak l}^v u-yz \pmod{\mathfrak l^n},
    \end{align*}
    so we let $\alpha=x-t/2$. We split into cases according to $v_{\mathfrak l}(y)=j=0,1,\ldots,n$, where $v_{\mathfrak l}(y)=n$ indicates $y=0$ in $\mathcal O_{E,\mathfrak l}/\mathfrak l^n \mathcal O_{E,\mathfrak l}$. For a given $y$ with $v_{\mathfrak l}(y)=j$, the equation has a solution if and only if
    \begin{align}\label{eqn_R_j condition}
        \pi_{\mathfrak l}^v u -\alpha^2 \in \mathfrak l^j.
    \end{align}
    Further, for each such $\alpha$, the equation is equivalent to
    \begin{align}\label{eqn_z equation}
        z \equiv (y')^{-1}\pi_{\mathfrak l}^{-j}(\pi_{\mathfrak l}^v u-\alpha^2) \pmod{\mathfrak l^{n-j}},
    \end{align}
    where $y=\pi_{\mathfrak l}^j y'$, $y' \in (\mathcal O_{E,\mathfrak l}/\mathfrak l^n \mathcal O_{E,\mathfrak l})^{\times}$.
    There are $q_{\mathfrak l}^j$ lifts $z$ modulo $\mathfrak l^n$ satisfying \eqref{eqn_z equation}. 
    Let $R_j$ be the number of $\alpha \in \mathcal O_{E,\mathfrak l}/\mathfrak l^n \mathcal O_{E,\mathfrak l}$ satisfying \eqref{eqn_R_j condition}. Note that the number of $y \in \mathcal O_{E,\mathfrak l}/\mathfrak l^n \mathcal O_{E,\mathfrak l}$ with $v_{\mathfrak l}(y)=j$ is
    \begin{align*}
        \begin{cases}
            (q_{\mathfrak l}-1)q_{\mathfrak l}^{n-j-1} & \text{if} \ j<n, \\
            1 & \text{if} \ j=n.
        \end{cases}
    \end{align*}
    Thus we have
    \begin{align*}
        |C_{\mathfrak l^n}(d,t)|=\sum_{j=0}^{n-1} (q_{\mathfrak l}-1)q_{\mathfrak l}^{n-j-1}R_j q_{\mathfrak l}^j + R_n q_{\mathfrak l}^n.
    \end{align*}

    Now we compute $R_j$. When $j \leq v$, the condition \eqref{eqn_R_j condition} is equivalent to $\alpha \equiv 0 \pmod{\mathfrak l^{\lceil j/2 \rceil}}$. So the number of such $\alpha \in \mathcal O_{E,\mathfrak l}/\mathfrak l^n \mathcal O_{E,\mathfrak l}$ is $q_{\mathfrak l}^{n-\lceil j/2 \rceil}$. When $j>v$, the condition \eqref{eqn_R_j condition} holds only if $v$ is even and $u$ is a square modulo $\mathfrak l$. In this case, we may write $\alpha=\pi_{\mathfrak l}^{v/2}\alpha'$, so the condition is equivalent to
    \begin{align*}
        \alpha'^2 \equiv u \pmod{\mathfrak l^{j-v}}.
    \end{align*}
    There are two $\alpha'$ modulo $\mathfrak l^{j-v}$ satisfying the above equation, and for each solution we have $\alpha=\pi_{\mathfrak l}^{v/2}\alpha'$ modulo $\mathfrak l^{j-v/2}$. Each $\alpha$ has $q_{\mathfrak l}^{n-j+v/2}$ lifts modulo $\mathfrak l^n$. In short,
    \begin{align*}
        R_j = \begin{cases}
            q_{\mathfrak l}^{n-\lceil j/2 \rceil} & \text{if} \ 0\leq j \leq v, \\
            2q_{\mathfrak l}^{n-j+v/2} & \text{if} \ j>v, \ \text{$v$ is even, and $u$ is a square modulo $\mathfrak l$}, \\
            0 & \text{otherwise}.
        \end{cases}
    \end{align*}
    Therefore, if $v<n$ so that the exceptional case $j=n$ lies in the range $v<j$, we have
    \begin{align*}
        |C_{\mathfrak l^n}(d,t)|=\sum_{j=0}^v (q_{\mathfrak l}-1)q_{\mathfrak l}^{2n-\lceil j/2 \rceil-1}+2\delta_{v,u}q_{\mathfrak l}^{n+v/2}+2\delta_{v,u}(q_{\mathfrak l}-1)q_{\mathfrak l}^{n-1}\sum_{j=v+1}^{n-1}q_{\mathfrak l}^{n-j+v/2},
    \end{align*}
    where $\delta_{v,u}=1$ if $v$ is even and $u$ is a square modulo $\mathfrak l$, and $\delta_{v,u}=0$ otherwise. If $v=n$, so that no $j$ satisfies $j>v$, then 
    \begin{align*}
        |C_{\mathfrak l^n}(d,t)|=\sum_{j=0}^{n-1}(q_{\mathfrak l}-1)q_{\mathfrak l}^{2n-\lceil j/2 \rceil-1} +  q_{\mathfrak l}^{n+\lfloor n/2 \rfloor}.
    \end{align*}
    One can verify that
    \begin{align*}
        2(q_{\mathfrak l}-1)q_{\mathfrak l}^{n-1}\sum_{j=v+1}^{n-1}q_{\mathfrak l}^{n-j+v/2}=2(q_{\mathfrak l}^{2n-v/2-1}-q_{\mathfrak l}^{n+v/2}),
    \end{align*}
    which completes the proof for $|C_{\mathfrak l^n}(d,t)|$ for $\ell>2$.

    For $\ell=2$, we use the following simpler argument. Consider the equation
    \begin{align}\label{eqn_x^2-tx+d=0}
        x^2-tx+d \equiv 0 \pmod{\mathfrak l^j}.
    \end{align}
    Let $P(x)=x^2-tx+d$. If $x_0 \in \mathcal O_{E,\mathfrak l}/\mathfrak l^j \mathcal O_{E,\mathfrak l}$ is a root of \eqref{eqn_x^2-tx+d=0}, then for any root $x \in \mathcal O_{E,\mathfrak l}/\mathfrak l^j \mathcal O_{E,\mathfrak l}$, we have
    \begin{align*}
        P(x)-P(x_0)=(x-x_0)(x+x_0-t) \equiv 0 \pmod{\mathfrak l^j},
    \end{align*}
    so either
    \begin{align*}
        x \equiv x_0 \pmod{\mathfrak l^{\lceil j/2 \rceil}}
    \end{align*}
    or
    \begin{align*}
        x \equiv -x_0+t \pmod{\mathfrak l^{\lceil j/2 \rceil}}.
    \end{align*}
    Hence the number of solutions of \eqref{eqn_x^2-tx+d=0} is $\leq 2 q_{\mathfrak l}^{\lfloor j/2 \rfloor}$. Recall that we want to count the number of triples $(x,y,z) \in (\mathcal O_{E,\mathfrak l}/\mathfrak l^n \mathcal O_{E,\mathfrak l})^3$ with
    \begin{align*}
        x^2-tx+d \equiv yz \pmod{\mathfrak l^n}.
    \end{align*}
    If we let $v_{\mathfrak l}(y)=j$, then the root $x$ must satisfy \eqref{eqn_x^2-tx+d=0}. The number of pairs $(y,z)$ with $v_{\mathfrak l}(y)=j$ is given by the same formula as in the case $\ell>2$, so we have
    \begin{align*}
        |C_{\mathfrak l^n}(d,t)| \leq \sum_{j=0}^{n-1} 2q_{\mathfrak l}^{\lfloor j/2 \rfloor}(q_{\mathfrak l}-1)q_{\mathfrak l}^{n-j-1}q_{\mathfrak l}^n + 2q_{\mathfrak l}^{\lfloor j/2 \rfloor}q_{\mathfrak l}^n,
    \end{align*}
    which is $O(q_{\mathfrak l}^{2n})$.

    For the remaining assertion, note that $A_{\mathfrak l^n}(d)$ is a fiber of a homomorphism
    \begin{align*}
        \det: \GL_2(\mathcal O_{E,\mathfrak l}/\mathfrak l^n \mathcal O_{E,\mathfrak l}) \longrightarrow (\mathcal O_{E,\mathfrak l}/\mathfrak l^n \mathcal O_{E,\mathfrak l})^{\times},
    \end{align*}
    so the number of elements of $A_{\mathfrak l^n}(d)$ is the same as 
    \begin{align*}
        |\SL_2(\mathcal O_{E,\mathfrak l}/\mathfrak l^n \mathcal O_{E,\mathfrak l})|=\frac{|\GL_2(\mathcal O_{E,\mathfrak l}/\mathfrak l^n \mathcal O_{E,\mathfrak l})|}{|(\mathcal O_{E,\mathfrak l}/\mathfrak l^n \mathcal O_{E,\mathfrak l})^{\times}|}.
    \end{align*}
    Since we have the canonical surjection
    \begin{align*}
        \GL_2(\mathcal O_{E,\mathfrak l}/\mathfrak l^n \mathcal O_{E,\mathfrak l}) \longrightarrow \GL_2(\F_{q_{\mathfrak l}})
    \end{align*}
    with kernel
    \begin{align*}
        I+{\rm Mat}_2(\mathfrak l \mathcal O_{E,\mathfrak l}/\mathfrak l^n \mathcal O_{E,\mathfrak l}),
    \end{align*}
    the number of elements of $\GL_2(\mathcal O_{E,\mathfrak l}/\mathfrak l^n \mathcal O_{E,\mathfrak l})$ is $q_{\mathfrak l}^{4(n-1)}$. Therefore,
    \begin{align*}
        |\SL_2(\mathcal O_{E,\mathfrak l}/\mathfrak l^n \mathcal O_{E,\mathfrak l})|=\frac{q_{\mathfrak l}^{4(n-1)}}{q_{\mathfrak l}^{n-1}(q_{\mathfrak l}-1)}
    \end{align*}
    which completes the proof.
\end{proof}

To simplify the notation, we let
\begin{align*}
    \mathcal D_{\mathfrak l}&:=(\Z_{\ell}^{\times})^{k-1} \subset \mathcal O_{E,\mathfrak l}, \\
    \mathcal A_{\mathfrak l}&:=\{g \in \GL_2(\mathcal O_{E,\mathfrak l}): \det g \in \mathcal D_{\mathfrak l}\}
\end{align*}
and let $\bar{\mathcal D}_{\mathfrak l^n}$ and $\bar{\mathcal A}_{\mathfrak l^n}$ denote the images of $\mathcal D_{\mathfrak l}$ and $\mathcal A_{\mathfrak l}$ under reduction modulo $\mathfrak l^n$, respectively.

\begin{lemma}\label{lem_lower bound of open trace condition}
    For an open subgroup $\mathcal U$ of $\mathcal A_{\mathfrak l}$, denote by $\bar{\mathcal U}_{\mathfrak l^n}$ the image of $\mathcal U$ under mod $\mathfrak l^n$ reduction. Suppose that for $t \in O_{E,\mathfrak l}$, there exists $g_0 \in \mathcal U$ for which $\tr(g_0)=t$. Then we have
    \begin{align*}
        \frac{\#\{g \in \bar{\mathcal U}_{\mathfrak l^n}: \tr g \equiv t \pmod{\mathfrak l^n}\}}{|\bar{\mathcal U}_{\mathfrak l^n}|} \gg_{\mathcal U,t} \frac{1}{q_{\mathfrak l}^n}.
    \end{align*}
\end{lemma}

\begin{proof}
    Note that we may assume $(g_0)_{21} \neq 0$. Indeed, if not, then we replace $g_0$ as follows. If $g_0 \neq zI$ for any $z$, we consider
    \begin{align*}
        \begin{pmatrix}
            1 & 0 \\ -s & 1
        \end{pmatrix}
        g_0
        \begin{pmatrix}
            1 & 0 \\ s & 1
        \end{pmatrix}
    \end{align*}
    for $s \in \mathfrak l^e  \mathcal O_{E,\mathfrak l}$. It belongs to $\mathcal U$, its trace is the same as $\tr(g_0)$, and its $(2,1)$-component is given by $s((g_0)_{22}-(g_0)_{11})-s^2 (g_0)_{12}$, so we may choose a suitable $s$ for which the $(2,1)$-component is nonzero. If $g_0=zI$ for some $z$, then we may consider
    \begin{align*}
        g_0\begin{pmatrix}
            1 & 0 \\ s & 1
        \end{pmatrix}
    \end{align*}
    for $s \in \mathfrak l^e  \mathcal O_{E,\mathfrak l}$ so that $g_0' \in \mathcal U$. The trace is given by $\tr(g_0)$, and $(g_0')_{21}=zs \neq 0$.
    
    Let $v_0:=v_{\mathfrak l}((g_0)_{21})$ be the $\mathfrak l$-adic valuation of $(g_0)_{21}$.
    
    We choose an open subgroup $\mathcal V$ of $\mathcal D_{\mathfrak l}$ that is contained in the intersection of $1+\mathfrak l^L \mathcal O_{E,\mathfrak l}$ and the preimage of $\mathcal U$ under the map
    \begin{align*}
        \iota_{\mathfrak l}: \mathcal D_{\mathfrak l} \longrightarrow \mathcal A_{\mathfrak l}
    \end{align*}
    defined by $\iota_{\mathfrak l}(d)=\begin{pmatrix}
        d & 0 \\ 0 & 1
    \end{pmatrix}.$ Here, $L>e$ is a sufficiently large positive integer depending on $\mathcal U$ and $t$, which will be specified below. Denote by $\bar{\mathcal V}_{\mathfrak l^n} \subset \bar{\mathcal D}_{\mathfrak l^n}$ the image of $\mathcal V$ under mod $\mathfrak l^n$ reduction. We will use $g_0$ to construct sufficiently many elements $g \in \bar{\mathcal U}_{\mathfrak l^n}$ satisfying $\tr g \equiv t \pmod{\mathfrak l^n}$.

    Let
    \begin{align*}
        l(x):=\begin{pmatrix}
            1 & 0 \\ x & 1
        \end{pmatrix},
    \qquad 
        u(y):=\begin{pmatrix}
            1 & y \\ 0 & 1
        \end{pmatrix},
    \qquad
        a(z):=\begin{pmatrix}
            1+z & 0 \\ 0 & (1+z)^{-1}
        \end{pmatrix}
    \end{align*}
    for $x,y,z \in \mathcal O_{E,\mathfrak l}$.
    We consider a matrix of the form
    \begin{align*}
        g=g_0 \iota_{\mathfrak l}(d) l(x) a(z) u(y) \in \GL_2(\mathcal O_{E,\mathfrak l}).
    \end{align*}
    Note that if $x,y,z \in \mathfrak l^e \mathcal O_{E,\mathfrak l}$, then $l(x),u(y),a(z) \in \mathcal U$. Thus, if further $d \in \mathcal V \subset 1+\mathfrak l^L \mathcal O_{E,\mathfrak l}$, then $g \in \mathcal U$.

    To simplify the notation, write $g_0(d,x,z)=g_0\iota_{\mathfrak l}(d)l(x)a(z)$. It is easy to verify that $d,x,z$ are uniquely determined by $g_0(d,x,z)$. We also note that, in the $\mathfrak l$-adic topology,
    \begin{align*}
        \iota_{\mathfrak l}(d)l(x)a(z) \longrightarrow I
    \end{align*}
    as
    \begin{align*}
        d \longrightarrow 1, \qquad x \longrightarrow 0, \qquad z \longrightarrow 0.
    \end{align*}
    Consequently, by enlarging $L$ if necessary, we have
    \begin{align*}
        \tr(g_0(d,x,z)-g_0)&=\tr(g_0(d,x,z))-t \in \mathfrak l^{v_0+e}\mathcal O_{E,\mathfrak l}, \\
        v_{\mathfrak l}(g_0(d,x,z)_{21})&=v_0
    \end{align*}
    if $x,z \in \mathfrak l^L \mathcal O_{E,\mathfrak l}$ and $d \in \mathcal V$.

    As
    \begin{align*}
        \tr(g_0(d,x,z)u(y))=\tr(g_0(d,x,z))+g_0(d,x,z)_{21}y,
    \end{align*}
    the two conditions
    \begin{align*}
        \tr(g) \equiv t \pmod{\mathfrak l^n}
    \end{align*}
    and
    \begin{align}\label{eqn_tr g_0 equiv t - g_0 y}
        \tr(g_0(d,x,z)) \equiv t - g_0(d,x,z)_{21}y \pmod{\mathfrak l^n}
    \end{align}
    are equivalent. If we write $g_0(d,x,z)_{21}=\pi_{\mathfrak l}^{v_0}g_0'$ for $g_0' \in \mathcal O_{E,\mathfrak l}^{\times}$, the equation \eqref{eqn_tr g_0 equiv t - g_0 y} can be rewritten as
    \begin{align*}
        g_0' y \equiv \pi_{\mathfrak l}^{-v_0}(t-\tr(g_0(d,x,z))) \pmod{\mathfrak l^{n-v_0}}. 
    \end{align*}
    Thus, for given $d \in \bar{\mathcal V}_{\mathfrak l^n}$ and $x,z \in \mathfrak l^L \mathcal O_{E,\mathfrak l}/ \mathfrak l^n \mathcal O_{E,\mathfrak l}$, there are $|\mathfrak l^{n-v_0} \mathcal O_{E,\mathfrak l}/ \mathfrak l^n \mathcal O_{E,\mathfrak l}|=q_{\mathfrak l}^{v_0}$ choices of $y$ mod $\mathfrak l^n$ that satisfy the desired equation. Further, these solutions $y$ satisfies $y \equiv 0 \pmod{\mathfrak l^e}$. Therefore, for $n \gg 1$, we have at least
    \begin{align*}
        |\bar{\mathcal V}_{\mathfrak l^n}|q_{\mathfrak l}^{2(n-L)+v_0}
    \end{align*}
    different elements $g(d,x,z,y)=g_0(d,x,z)u(y) \in \bar{\mathcal U}_{\mathfrak l^n}$ which satisfy $\tr(g) \equiv t \pmod{\mathfrak l^n}$.
    Note that, from the canonical surjection
    \begin{align*}
        \mathcal D_{\mathfrak l}/\mathcal V \longrightarrow \bar{\mathcal D}_{\mathfrak l^n}/\bar{\mathcal V}_{\mathfrak l^n},
    \end{align*}
    we have 
    \begin{align*}
        |\bar{\mathcal V}_{\mathfrak l^n}| \geq \frac{|\bar{\mathcal D}_{\mathfrak l^n}|}{[\mathcal D_{\mathfrak l}:\mathcal V]}.
    \end{align*}

    On the other hand, we have
    \begin{align*}
        |\bar{\mathcal U}_{\mathfrak l^n}| \leq |\bar{\mathcal A}_{\mathfrak l^n}|=|\bar{\mathcal D}_{\mathfrak l^n}||\SL_2(\mathcal O_{E,\mathfrak l}/\mathfrak l^n \mathcal O_{E,\mathfrak l})| = |\bar{\mathcal D}_{\mathfrak l^n}| q_{\mathfrak l}^{3n}(1-q_{\mathfrak l}^{-2}) \leq |\bar{\mathcal D}_{\mathfrak l^n}| q_{\mathfrak l}^{3n}.
    \end{align*}
    In conclusion, we obtain
    \begin{align*}
        \frac{\#\{g \in \bar{\mathcal U}_{\mathfrak l^n}: \tr(g) \equiv t \pmod{\mathfrak l^n}\}}{|\bar{\mathcal U}_{\mathfrak l^n}|} \geq \frac{1}{[\mathcal D_{\mathfrak l}:\mathcal V]} q_{\mathfrak l}^{-n-2L+v_0}.
    \end{align*}
    Since $[\mathcal D_{\mathfrak l}:\mathcal V]$, $L$, and $v_0$ depend on $\mathcal U$ and $t$, we get the assertion.
\end{proof}

For a rational prime $\ell$, let $\mathcal D_{\ell}$ be the diagonal image $\prod_{\mathfrak l \mid \ell} \iota_{\mathfrak l}((\Z_{\ell}^{\times})^{k-1})$ in $\prod_{\mathfrak l \mid \ell} \mathcal D_{\mathfrak l}$ and $\bar{\mathcal D}_{\prod_{\mathfrak l \mid \ell}\mathfrak l^{n_{\mathfrak l}}}$ be its mod $\prod_{\mathfrak l \mid \ell}\mathfrak l^{n_{\mathfrak l}}$ reduction. For $d \in \bar{\mathcal D}_{\prod_{\mathfrak l \mid \ell}\mathfrak l^{n_{\mathfrak l}}} \subset \prod_{\mathfrak l \mid \ell}(\mathcal O_{E,\mathfrak l}/\mathfrak l^{n_{\mathfrak l}} \mathcal O_{E,\mathfrak l})^{\times}$ and $t \in \prod_{\mathfrak l \mid \ell}\mathcal O_{E,\mathfrak l}/\mathfrak l^{n_{\mathfrak l}} \mathcal O_{E,\mathfrak l}$, we naturally define
    \begin{align*}
        C_{\prod_{\mathfrak l \mid \ell}\mathfrak l^{n_{\mathfrak l}}}(d,t)&=\{g \in \GL_2(\prod_{\mathfrak l \mid \ell}\mathcal O_{E,\mathfrak l}/\mathfrak l^{n_{\mathfrak l}} \mathcal O_{E,\mathfrak l}): \det g=d, \ \tr g = t\}, \\
        A_{\prod_{\mathfrak l \mid \ell}\mathfrak l^{n_{\mathfrak l}}}(d)&=\{g \in \GL_2(\prod_{\mathfrak l \mid \ell}\mathcal O_{E,\mathfrak l}/\mathfrak l^{n_{\mathfrak l}} \mathcal O_{E,\mathfrak l}): \det g=d\}.
    \end{align*}
\begin{lemma}\label{lem_total tr/det proportion}
    With the notation above, we have
    \begin{align*}
        &|C_{\prod_{\mathfrak l \mid \ell}\mathfrak l^{n_{\mathfrak l}}}(d,t)|=\prod_{\mathfrak l \mid \ell} |C_{\mathfrak l^{n_{\mathfrak l}}}(d,t_{\mathfrak l})| \ll \prod_{\mathfrak l \mid \ell} q_{\mathfrak l}^{2n_{\mathfrak l}}, \\
        &|A_{\prod_{\mathfrak l \mid \ell}\mathfrak l^{n_{\mathfrak l}}}(d,t)|=\prod_{\mathfrak l \mid \ell} |A_{\mathfrak l^{n_{\mathfrak l}}}(d)| \gg \prod_{\mathfrak l \mid \ell} q_{\mathfrak l}^{3n_{\mathfrak l}-2}.
    \end{align*}
    Consequently,
    \begin{align*}
        &\#\{g \in \GL_2(\prod_{\mathfrak l \mid \ell}\mathcal O_{E,\mathfrak l}/\mathfrak l^{n_{\mathfrak l}} \mathcal O_{E,\mathfrak l}): \det g \in \bar{\mathcal D}_{\prod_{\mathfrak l \mid \ell}\mathfrak l^{n_{\mathfrak l}}}, \ \tr g = t\} \ll_{\ell} |\bar{\mathcal D}_{\prod_{\mathfrak l \mid \ell}\mathfrak l^{n_{\mathfrak l}}}|\prod_{\mathfrak l \mid \ell}q_{\mathfrak l}^{2n_{\mathfrak l}},   \\
        &|\bar{\mathcal A}_{\prod_{\mathfrak l \mid \ell}\mathfrak l^{n_{\mathfrak l}}}|=\#\{g \in \GL_2(\prod_{\mathfrak l \mid \ell}\mathcal O_{E,\mathfrak l}/\mathfrak l^{n_{\mathfrak l}} \mathcal O_{E,\mathfrak l}): \det g \in \bar{\mathcal D}_{\prod_{\mathfrak l \mid \ell}\mathfrak l^{n_{\mathfrak l}}}\} \gg_{\ell,E} |\bar{\mathcal D}_{\prod_{\mathfrak l \mid \ell}\mathfrak l^{n_{\mathfrak l}}}|\prod_{\mathfrak l \mid \ell} q_{\mathfrak l}^{3n_{\mathfrak l}}
    \end{align*}
    and hence
    \begin{align*}
        \frac{\#\{g \in \GL_2(\prod_{\mathfrak l \mid \ell}\mathcal O_{E,\mathfrak l}/\mathfrak l^{n_{\mathfrak l}} \mathcal O_{E,\mathfrak l}): \det g \in \bar{\mathcal D}_{\prod_{\mathfrak l \mid \ell}\mathfrak l^{n_{\mathfrak l}}}, \ \tr g = t\}}{|\bar{\mathcal A}_{\prod_{\mathfrak l \mid \ell}\mathfrak l^{n_{\mathfrak l}}}|} \ll_{\ell,E} \prod_{\mathfrak l \mid \ell} \frac{1}{q_{\mathfrak l}^{n_{\mathfrak l}}}.
    \end{align*}
\end{lemma}

\begin{proof}
    This is an immediate consequence of Lemma~\ref{lem_A lambda^n d, t number} for $\ell>2$. 
\end{proof}

If $\mathfrak l \gg_F 1$, then one has
\begin{align*}
    F(t) \equiv 0 \pmod{\mathfrak l^n} \qquad \text{if and only if} \qquad t\equiv t_i \pmod{\mathfrak l^{\lceil n/e_i \rceil}}
\end{align*}
for some root $t_i$ of $F_i$ in $\bar E$. In fact, since we are considering roots modulo powers of $\mathfrak l$, the situation is clearer if we consider the roots in $\mathcal O_{E,\mathfrak l}$. As we assumed $\ell \gg_F 1$, each root of $F_i$ mod $\mathfrak l$ is simple. The condition $F(t) \equiv 0 \pmod{\mathfrak l^n}$ implies $F_i(\bar t)=0$ for some $i$, where $\bar t$ is the mod $\mathfrak l$ reduction of $t$, and hence by Hensel's lemma, there exists a root $t_{\mathfrak l} \in \mathcal O_{E,\mathfrak l}$ of $F$ with $t_{\mathfrak l,j} \equiv \bar t \pmod{\mathfrak l}$. Under the natural embedding $\bar E \hookrightarrow \bar E_{\mathfrak l}$, each of these roots $t_{\mathfrak l}$ is one of the $t_i \in \bar E \cap \mathcal O_{E,\mathfrak l}$. In conclusion, if we let $t_{\mathfrak l,1},\ldots,t_{\mathfrak l,r_{\mathfrak l}} \in \mathcal O_{E,\mathfrak l}$ be the distinct roots of $F$, and $i_j$ be the index for which $t_{\mathfrak l,j}$ is a root of $F_{i_j}$, then
\begin{align*}
    F(t) \equiv 0 \pmod{\mathfrak l^n} \qquad \text{if and only if} \qquad t\equiv t_{\mathfrak l,j} \pmod{\mathfrak l^{\lceil n/e_{i_j} \rceil}}
\end{align*}
for exactly one $t_{\mathfrak l,j}$ among $t_{\mathfrak l,1},\ldots, t_{\mathfrak l,r_{\mathfrak l}}$. An implicit way to express $e_{i_j}$ is as follows. The number $e_{i_j}$ equals $m_{\mathfrak l,j}$, the exponent for which $F$ can be written as
\begin{align*}
    F(T)=(T-t_{\mathfrak l,j})^{m_{\mathfrak l,j}}H_{\mathfrak l,j}(T),
\end{align*}
with $H_{\mathfrak l,j}(t_{\mathfrak l,j}) \neq 0$. Since the roots of $F_i$ and $F_j$ do not overlap, we have $H_{\mathfrak l,j}(t_{\mathfrak l,j}) \in \mathcal O_{E,\mathfrak l}^{\times}$.

However, this does not hold in general, as the roots of mod $\mathfrak l^n$ reductions of $F_i$ and $F_j$ can coincide when $\mathfrak l$ is not sufficiently large. The following lemma provides the precise statement for the general situation.

\begin{lemma}\label{lem_decomposing the condition F(tr g)=0}
    Let $\mathfrak l$ be a fixed prime of $E$. Let $t_{\mathfrak l,1}, \ldots, t_{\mathfrak l,r_{\mathfrak l}} \in \mathcal O_{E,\mathfrak l}$ be the distinct roots of $F$, and for each $1 \leq j\leq r_{\mathfrak l}$, write
    \begin{align*}
        F(T)=(T-t_{\mathfrak l,j})^{m_{\mathfrak l,j}}H_{\mathfrak l,j}(T)
    \end{align*}
    with $H_{\mathfrak l,j}(t_{\mathfrak l,j}) \neq 0$. Let
    \begin{align*}
        c_{\mathfrak l,j}=v_{\mathfrak l}(H_{\mathfrak l,j}(t_{\mathfrak l,j})).
    \end{align*}
    There exists $N_{\mathfrak l}>0$ such that if $n>N_{\mathfrak l}$, then
    \begin{align*}
        \{t \in \mathcal O_{E,\mathfrak l}: F(t) \equiv 0 \pmod{\mathfrak l^n}\} = \bigsqcup_{1 \leq j \leq r_{\mathfrak l}}\{t \in \mathcal O_{E,\mathfrak l}: t \equiv t_{\mathfrak l,j} \pmod{\mathfrak l^{b_{\mathfrak l,j}(n)}}\},
    \end{align*}
    where 
    \begin{align*}
        b_{\mathfrak l,j}(n):=\bigg\lceil \frac{n-c_{\mathfrak l,j}}{m_{\mathfrak l,j}}\bigg\rceil.    
    \end{align*}
\end{lemma}

\begin{proof}
    By choosing a suitably large $L_{\mathfrak l}>0$, we may take compact open balls
    \begin{align*}
        t_{\mathfrak l,j}+\mathfrak l^{L_{\mathfrak l}}\mathcal O_{E,\mathfrak l}
    \end{align*}
    for $1 \leq j \leq r_{\mathfrak l}$ which are pairwise disjoint. Since $H_{\mathfrak l,j}$ is continuous, we may assume that
    \begin{align*}
        v_{\mathfrak l}(H_{\mathfrak l,j}(t))=v_{\mathfrak l}(H_{\mathfrak l,j}(t_{\mathfrak l,j}))=c_{\mathfrak l,j}
    \end{align*}
    if $t \in t_{\mathfrak l,j}+\mathfrak l^{L_{\mathfrak l}}\mathcal O_{E,\mathfrak l}$. Note that the complement
    \begin{align*}
        \mathcal O_{E,\mathfrak l} \setminus \bigsqcup_{j=1}^{r_{\mathfrak l}}(t_{\mathfrak l,j}+\mathfrak l^{L_{\mathfrak l}}\mathcal O_{E,\mathfrak l})
    \end{align*}
    is compact and $F(t) \neq 0$ for any $t$ in this complement obviously, so the value $|F(t)|_{\mathfrak l}$ is bounded below by some positive number. In other words, $v_{\mathfrak l}(F(t))$ is bounded above for $t \not\in \bigsqcup_{j=1}^{r_{\mathfrak l}}(t_{\mathfrak l,j}+\mathfrak l^{L_{\mathfrak l}}\mathcal O_{E,\mathfrak l})$. Thus for sufficiently large $N_{\mathfrak l}$, if $n>N_{\mathfrak l}$, then $F(t) \equiv 0 \pmod{\mathfrak l^n}$ implies that $t$ belongs to $t_{\mathfrak l,j}+\mathfrak l^{L_{\mathfrak l}}\mathcal O_{E,\mathfrak l}$ for some $j.$ Consequently,
    \begin{align*}
        n \leq v_{\mathfrak l}(F(t))=m_{\mathfrak l,j}v_{\mathfrak l}(t-t_{\mathfrak l,j})+v_{\mathfrak l}(H_{\mathfrak l,j}(t))=m_{\mathfrak l,j}v_{\mathfrak l}(t-t_{\mathfrak l,j})+c_{\mathfrak l,j},
    \end{align*}
    and hence 
    \begin{align*}
        v_{\mathfrak l}(t-t_{\mathfrak l,j}) \geq b_{\mathfrak l,j}(n).
    \end{align*}
    This is equivalent to saying that $t$ satisfies $t\equiv t_{\mathfrak l,j} \pmod{\mathfrak l^{b_{\mathfrak l,j}(n)}}$.

    The converse is easy to verify; if $t \equiv t_{\mathfrak l,j} \pmod{\mathfrak l^{b_{\mathfrak l,j}(n)}}$ for some $1 \leq j \leq r_{\mathfrak l}$, then $F(t)=(t-t_{\mathfrak l,j})^{m_{\mathfrak l,j}}H_{\mathfrak l,j}(t)$ is divisible by $\pi_{\mathfrak l}^{m_{\mathfrak l,j}b_{\mathfrak l,j}(n)+c_{\mathfrak l,j}}$, and it is obvious that the exponent $m_{\mathfrak l,j}b_{\mathfrak l,j}(n)+c_{\mathfrak l,j}$ is larger than or equal to $n$. 
\end{proof}

Let $\mathcal A_{\ell}$ be the subset of $\prod_{\mathfrak l \mid \ell} \mathcal A_{\mathfrak l}$ consisting of the elements whose determinant belongs to $\mathcal D_{\ell}$.
Let $S$ be a finite set of rational primes, set $\mathcal A_S:=\prod_{\ell \in S} \mathcal A_{\ell}$, and let $\mathcal U \subset \mathcal A_{S}$ be an open subgroup. For $\mathfrak o$ such that every rational prime $\ell$ dividing $\mathfrak o$ belongs to $S$, denote by $\bar{\mathcal U}_{\mathfrak o}$ the mod $\mathfrak o$ reduction of $\mathcal U$.
Let $\mathcal U_{\mathfrak l}$ be the $\mathfrak l$-adic projective image of $\mathcal U$ and
\begin{align*}
    J_{\mathfrak l}=\{1 \leq j \leq r_{\mathfrak l}: \text{there exists $g_0 \in \mathcal U_{\mathfrak l}$ such that $\tr (g_0)=t_{\mathfrak l,j}$}\}.
\end{align*}
If $j \not\in J_{\mathfrak l}$, then there is an integer $M_{\mathfrak l,j}$ such that
\begin{align*}
    \tr (g) \not\equiv t_{\mathfrak l,j} \pmod{\mathfrak l^{M_{\mathfrak l,j}}}
\end{align*}
for every $g \in \mathcal U_{\mathfrak l}$. Indeed, if not, then we may choose $g_i$ for $i \geq 1$ for which $\tr (g_i) \to t_{\mathfrak l,j}$ in $\mathfrak l$-adic topology, which implies $t_{\mathfrak l,j} \in \tr(\mathcal U_{\mathfrak l})$ as $\tr(\mathcal U_{\mathfrak l})$ is closed. Thus, if one chooses $N_{\mathfrak l}$ for which $b_{\mathfrak l,j}(N_{\mathfrak l}$ is larger than $M_{\mathfrak l,j}$ for all $1 \leq j \leq r_{\mathfrak l}$, then the existence of $g_0$ with $\tr(g_0) \equiv t_{\mathfrak l,j} \pmod{\mathfrak l^{b_{\mathfrak l,j}(N_{\mathfrak l})}}$ implies $j \in J_{\mathfrak l}$. 

We address the `bad prime part' of the essential multiplicativity in the following proposition. Let
    \begin{align*}
        \mathcal C_{F,\mathfrak o}(\mathcal U):=\{ g \in \bar{\mathcal U}_{\mathfrak o} : F(\tr (g)) \equiv 0 \pmod{\mathfrak o}\}
    \end{align*}
    and
    \begin{align*}
        h_{F}(\mathfrak o; \mathcal U):=\frac{|\mathcal C_{F,\mathfrak o}(\mathcal U)|}{|\bar{\mathcal U}_{\mathfrak o}|}.
    \end{align*}

\begin{proposition}\label{prop_local density mod forms bad prime essential multiplicativity}
    Suppose that an integral ideal $\mathfrak o$ is given by
    \begin{align*}
        \mathfrak o=\prod_{\mathfrak l \mid \mathfrak o} \mathfrak l^{v_{\mathfrak l}}.
    \end{align*}
    Assume that if a rational prime $\ell$ divides $\mathfrak o$, then $\ell \in S$.
    We have
    \begin{align*}
        h_F(\mathfrak o; \mathcal U) \ll_{F,\mathcal U} \prod_{\mathfrak l \mid \mathfrak o} h_{F}(\mathfrak l^{v_{\mathfrak l}};\mathcal U).
    \end{align*}
    The implied constant is independent of the exponents $v_{\mathfrak l}$.
\end{proposition}

\begin{proof}

    Note that if $h_F(\mathfrak o;\mathcal U)=0$, the assertion holds obviously. We assume $h_F(\mathfrak o;\mathcal U) \neq 0$ and hence $h_F(\mathfrak l^{v_{\mathfrak l}};\mathcal U) \neq 0$ for all $\mathfrak l \mid \mathfrak o$.

    Let $N_0$ be a positive integer which is larger than $N_{\mathfrak l}$ for all $\mathfrak l \mid \mathfrak o$, where $N_{\mathfrak l}$ is chosen by Lemma~\ref{lem_decomposing the condition F(tr g)=0}. Note that for $\mathfrak o' \mid \mathfrak o$, if we write the mod $\mathfrak o'$ reduction as
    \begin{align*}
        {\rm red}_{\mathfrak o, \mathfrak o'}: \mathcal U_{\mathfrak o} \longrightarrow \mathcal U_{\mathfrak o'},
    \end{align*}
    then we have $\mathcal C_{F,\mathfrak o}(\mathcal U) \subseteq {\rm red}_{\mathfrak o, \mathfrak o'}^{-1}(\mathcal C_{F,\mathfrak o}(\mathcal U))$ so that
    \begin{align*}
        h_F(\mathfrak o;\mathcal U) = \frac{|\mathcal C_{F,\mathfrak o}(\mathcal U)|}{|\bar{\mathcal U}_{\mathfrak o}|} \leq \frac{|{\rm red}_{\mathfrak o, \mathfrak o'}^{-1}(\mathcal C_{F,\mathfrak o}(\mathcal U))|}{|\bar{\mathcal U}_{\mathfrak o}|} = \frac{|\mathcal C_{F,\mathfrak o}(\mathcal U)|}{|\bar{\mathcal U}_{\mathfrak o'}|} = h_F(\mathfrak o';\mathcal U).
    \end{align*}
    In particular,
    \begin{align*}
        h_F(\mathfrak o;\mathcal U) \leq h_F \bigl(\prod_{\substack{\mathfrak l \mid \mathfrak o \\ v_{\mathfrak l} > N_0}}\mathfrak l^{v_{\mathfrak l}};\mathcal U\bigl).
    \end{align*}
    Since there are only finitely many $\ell$ that belong to $S$, we also have
    \begin{align*}
        1 \ll_{F,S} \prod_{\substack{\mathfrak l \mid \mathfrak o \\ v_{\mathfrak l} \leq  N_0}}h_F(\mathfrak l^{v_{\mathfrak l}};\mathcal U),
    \end{align*}
    so it suffices to prove the assertion for $\mathfrak o$ such that every $v_{\mathfrak l}$ is larger than $N_0$.

    By Lemma~\ref{lem_decomposing the condition F(tr g)=0}, if $g=(g_{\mathfrak l})_{\mathfrak l \mid \mathfrak o} \in \bar{\mathcal U}_{\mathfrak o}$ satisfies $F(\tr g_{\mathfrak l}) \equiv 0 \pmod{\mathfrak l^{v_{\mathfrak l}}}$, then $\tr g_{\mathfrak l} \equiv t_{\mathfrak l,j} \pmod{\mathfrak l^{b_{\mathfrak l,j}(v_{\mathfrak l})}}$ for some $j_{\mathfrak l} \in J_{\mathfrak l}.$ This implies that
    \begin{align*}
        h_F(\mathfrak o;\mathcal U) \leq \sum_{(j_{\mathfrak l}) \in \prod_{\mathfrak l \mid \mathfrak o} J_{\mathfrak l}} \frac{\#\{g \in \bar{\mathcal U}_{\mathfrak o}: \tr g_{\mathfrak l} \equiv t_{\mathfrak l,j_{\mathfrak l}} \pmod{\mathfrak l^{b_{\mathfrak l,j}(v_{\mathfrak l})}}\}}{|\bar{\mathcal U}_{\mathfrak o}|}.
    \end{align*}
   To simplify the notation, set $b_{\mathfrak l}=b_{\mathfrak l,j_{\mathfrak l}}(v_{\mathfrak l})$ and $\mathfrak b:=(\mathfrak l^{b_{\mathfrak l}})_{\mathfrak l \mid \mathfrak o}$. As $b_{\mathfrak l} \leq v_{\mathfrak l}$, we have
    \begin{align*}
        \frac{\#\{g \in \bar{\mathcal U}_{\mathfrak o}: \tr g_{\mathfrak l} \equiv t_{\mathfrak l,j_{\mathfrak l}} \pmod{\mathfrak l^{b_{\mathfrak l}}}\}}{|\bar{\mathcal U}_{\mathfrak o}|} \leq \frac{\#\{g \in \bar{\mathcal U}_{\mathfrak b}: \tr g_{\mathfrak l} \equiv t_{\mathfrak l,j_{\mathfrak l}} \pmod{\mathfrak l^{b_{\mathfrak l}}}\}}{|\bar{\mathcal U}_{\mathfrak b}|}
    \end{align*}
    Note that
    \begin{align*}
        \frac{\#\{g \in \bar{\mathcal U}_{\mathfrak b}: \tr g_{\mathfrak l} \equiv t_{\mathfrak l,j_{\mathfrak l}} \pmod{\mathfrak l^{b_{\mathfrak l}}}\}}{|\bar{\mathcal U}_{\mathfrak b}|} &\leq [\bar{\mathcal A}_{\mathfrak b}:\bar{\mathcal U}_{\mathfrak b}] \frac{\#\{g \in \bar{\mathcal A}_{\mathfrak b}: \tr g_{\mathfrak l} \equiv t_{\mathfrak l,j_{\mathfrak l}} \pmod{\mathfrak l^{b_{\mathfrak l}}}\}}{|\bar{\mathcal A}_{\mathfrak b}|} \\
        & \leq [\mathcal A_S:\mathcal U] \frac{\#\{g \in \bar{\mathcal A}_{\mathfrak b}: \tr g_{\mathfrak l} \equiv t_{\mathfrak l,j_{\mathfrak l}} \pmod{\mathfrak l^{b_{\mathfrak l}}}\}}{|\bar{\mathcal A}_{\mathfrak b}|},
    \end{align*}
    where the last inequality follows again from the canonical surjection
    \begin{align*}
        \mathcal A_S/\mathcal U \longrightarrow \mathcal A_{\mathfrak b}/\mathcal U_{\mathfrak b}.
    \end{align*}
    By the definition $\mathcal A_S=\prod_{\ell \in S} \mathcal A_{\ell}$ and the Chinese remainder theorem, 
    \begin{align*}
        \frac{\#\{g \in \bar{\mathcal A}_{\mathfrak b}: \tr g_{\mathfrak l} \equiv t_{\mathfrak l,j_{\mathfrak l}} \pmod{\mathfrak l^{b_{\mathfrak l}}}\}}{|\bar{\mathcal A}_{\mathfrak b}|} = \prod_{\ell \in S} \frac{\#\{g \in \bar{\mathcal A}_{\mathfrak b_{\ell}}: \tr g_{\mathfrak l} \equiv t_{\mathfrak l,j_{\mathfrak l}} \pmod{\mathfrak l^{b_{\mathfrak l}}}\}}{|\bar{\mathcal A}_{\mathfrak b_{\ell}}|},
    \end{align*}
    where $\mathfrak b_{\ell}:=\prod_{\mathfrak l \mid \ell} \mathfrak l^{b_{\mathfrak l}}$.
    According to Lemmas~\ref{lem_lower bound of open trace condition} and ~\ref{lem_total tr/det proportion},
    \begin{align*}
        \frac{\#\{g \in \bar{\mathcal U}_{\mathfrak l^{b_{\mathfrak l}}}: \tr g_{\mathfrak l} \equiv t_{\mathfrak l,j_{\mathfrak l}} \pmod{\mathfrak l^{b_{\mathfrak l}}}\}}{|\bar{\mathcal A}_{\mathfrak b_{\ell}}|} &\ll_{S,E} \prod_{\mathfrak l \mid \ell} q^{-b_{\mathfrak l}} \\
        &\ll_{F,\mathcal U} 
        \prod_{\mathfrak l \mid \ell} \frac{\#\{g \in \bar{\mathcal U}_{\mathfrak l^{b_{\mathfrak l}}}: \tr g \equiv t_{\mathfrak l,j_{\mathfrak l}} \pmod{\mathfrak l^{b_{\mathfrak l}}}\}}{|\bar{\mathcal U}_{\mathfrak l^{b_{\mathfrak l}}}|}
    \end{align*}
    and consequently,
    \begin{align*}
        h_F(\mathfrak o;\mathcal U) \ll_{F,\mathcal U} \sum_{(j_{\mathfrak l}) \in \prod_{\mathfrak l \mid \mathfrak o} J_{\mathfrak l}} \prod_{\mathfrak l \mid \ell} \frac{\#\{g \in \bar{\mathcal U}_{\mathfrak l^{b_{\mathfrak l}}}: \tr g \equiv t_{\mathfrak l,j_{\mathfrak l}} \pmod{\mathfrak l^{b_{\mathfrak l}}}\}}{|\bar{\mathcal U}_{\mathfrak l^{b_{\mathfrak l}}}|}.
    \end{align*}
    For any $\mathfrak l$, by Lemma~\ref{lem_decomposing the condition F(tr g)=0} again,
    \begin{align*}
        \frac{\#\{g \in \bar{\mathcal U}_{\mathfrak l^{b_{\mathfrak l}}}: \tr g \equiv t_{\mathfrak l,j_{\mathfrak l}} \pmod{\mathfrak l^{b_{\mathfrak l}}}\}}{|\bar{\mathcal U}_{\mathfrak l^{b_{\mathfrak l}}}|} 
        &=\frac{\#\{g \in \bar{\mathcal U}_{\mathfrak l^{v_{\mathfrak l}}}: \tr g \equiv t_{\mathfrak l,j_{\mathfrak l}} \pmod{\mathfrak l^{b_{\mathfrak l}}}\}}{|\bar{\mathcal U}_{\mathfrak l^{v_{\mathfrak l}}}|}\\
        &\leq h_F(\mathfrak l^{v_{\mathfrak l}}; \mathcal U).
    \end{align*}
    We complete the proof by noting that the number of elements of $\prod_{\mathfrak l \mid \mathfrak o} J_{\mathfrak l}$ is $O_{F,S}(1)$.
\end{proof}

We now consider a similar property for a single large rational prime $\ell$ as follows. The estimate is stronger than that in Proposition~\ref{prop_local density mod forms bad prime essential multiplicativity}.

\begin{proposition}\label{prop_local density mod forms one good prime decomposition, stronger essetinal multiplicativity}
    Suppose $\ell \gg_{F,f} 1$. Let 
    \begin{align*}
        \mathfrak o = \prod_{\mathfrak l \mid \ell} \mathfrak l^{v_{\mathfrak l}}
    \end{align*}
    be an integral ideal of $E$ such that every prime factor lies over $\ell$. Then,
    \begin{align*}
        h_{F,f}(\mathfrak o)=(1+O(\ell^{-2}))\prod_{\mathfrak l \mid \ell} h_{F,f}(\mathfrak l^{v_{\mathfrak l}}).
    \end{align*}
\end{proposition}

\begin{proof}
    For $d \in \bar{\mathcal D}_{\mathfrak l^n}$, let
    \begin{align*}
        h_F(\mathfrak l^n; \mathcal A_{\mathfrak l}(d)):=\frac{\#\{g \in A_{\mathfrak l^n}(d): F(\tr g)\equiv 0 \pmod{\mathfrak l^n}\}}{|A_{\mathfrak l^n}(d)|}.
    \end{align*}
    Note that, by \cite[Theorem~2.2.2]{Loe17}, we have $h_{F,f}(\mathfrak l^n)=h_F(\mathfrak l^n;\mathcal A_{\mathfrak l})$. 
    As $|A_{\mathfrak l^n}(d)|=|A_{\mathfrak l^n}(1)|=:M_{\mathfrak l^n}$, we have
    \begin{align*}
        h_{F,f}(\mathfrak l^n)&=\frac{\sum_{d \in \bar{\mathcal D}_{\mathfrak l^n}} \#\{g \in A_{\mathfrak l^n}(d): F(\tr g)\equiv 0 \pmod{\mathfrak l^n}\}}{\sum_{d \in \bar{\mathcal D}_{\mathfrak l^n}} |A_{\mathfrak l^n}(d)|} \\
        &=\frac{\sum_{d \in \bar{\mathcal D}_{\mathfrak l^n}} h_F(\mathfrak l^n; \mathcal A_{\mathfrak l}(d)) M_{\mathfrak l^n}}{|\bar{\mathcal D}_{\mathfrak l^n}|M_{\mathfrak l^n}}
        = \frac{1}{|\bar{\mathcal D}_{\mathfrak l^n}|} \sum_{d \in \bar{\mathcal D}_{\mathfrak l^n}} h_F(\mathfrak l^n; \mathcal A_{\mathfrak l}(d)).
    \end{align*}
    In short, $h_{F,f}(\mathfrak l^n)$ is the average of $h_F(\mathfrak l^n; \mathcal A_{\mathfrak l}(d))$ over $d \in \bar{\mathcal D}_{\mathfrak l^n}$. 
    
    Note that
    \begin{align*}
        \{g \in A_{\mathfrak l^n}(d): F(\tr g)\equiv 0 \pmod{\mathfrak l^n}\} = \bigsqcup_{j=1}^{r_{\mathfrak l}} {\rm red}_{\mathfrak l^n,\mathfrak l^{b_{\mathfrak l,j}(n)}}^{-1}(C_{\mathfrak l^{b_{\mathfrak l,j}(n)}}(d,t_{\mathfrak l,j})).
    \end{align*}
    Since $\ell \gg_{F,f} 1$, $b_{\mathfrak l,j}(n)=\lceil n/m_{\mathfrak l,j}\rceil$. According to Lemma~\ref{lem_A lambda^n d, t number}, we obtain
    \begin{align*}
        h_F(\mathfrak l^n; \mathcal A_{\mathfrak l}(d))=\sum_{j=1}^{r_{\mathfrak l}} \frac{|C_{\mathfrak l^{b_{\mathfrak l,j}(n)}}(d)|}{|A_{\mathfrak l^{b_{\mathfrak l,j}(n)}}(d)|}=(1+O(q_{\mathfrak l}^{-1}))\sum_{j=1}^{r_{\mathfrak l}}q_{\mathfrak l}^{-b_{\mathfrak l,j}(n)}
    \end{align*}
    and consequently,
    \begin{align*}
        h_{F,f}(\mathfrak l^n)=(1+O(q_{\mathfrak l}^{-1}))\sum_{j=1}^{r_{\mathfrak l}}q_{\mathfrak l}^{-b_{\mathfrak l,j}(n)}.
    \end{align*}

    Suppose $h_{F,f}(\mathfrak l^n) \neq 0$. Put $\varepsilon_{\mathfrak l^n}(d):=(h_F(\mathfrak l^n; \mathcal A_{\mathfrak l}(d))-h_{F,f}(\mathfrak l^n))/h_{F,f}(\mathfrak l^n)$ so that $h_F(\mathfrak l^n; \mathcal A_{\mathfrak l}(d))=h_{F,f}(\mathfrak l^n)(1+\varepsilon_{\mathfrak l^n}(d))$. One can easily verify that
    \begin{align*}
        \varepsilon_{\mathfrak l^n}(d) \ll q_{\mathfrak l}^{-1}.
    \end{align*}
    Moreover, we have
    \begin{align*}
        |\bar{\mathcal D}_{\mathfrak l^n}|h_{F,f}(\mathfrak l^n)-\sum_{d \in \bar{\mathcal D}_{\mathfrak l^n}} h_F(\mathfrak l^n; \mathcal A_{\mathfrak l}(d)) = h_{F,f}(\mathfrak l^n) \sum_{d \in \bar{\mathcal D}_{\mathfrak l^n}}  \varepsilon_{\mathfrak l^n}(d) =0,
    \end{align*}

    Now consider 
    \begin{align*}
        h_{F,f}(\mathfrak o)=\frac{|\mathcal C_{F,f,\mathfrak o}|}{|\bar{\mathcal A}_{f,\mathfrak o}|}=\frac{\sum_{d \in \bar{\mathcal D}_{\mathfrak o}} \#\{g \in A_{\mathfrak o}(d): F(\tr g) \equiv 0 \pmod{\mathfrak o}\}}{\sum_{d \in \bar{\mathcal D}_{\mathfrak o}} |A_{\mathfrak o}(d)|}.
    \end{align*}
    Write $d=(d_{\mathfrak l})_{\mathfrak l \mid \ell} \in \bar{\mathcal D}_{\mathfrak o} \subset \prod_{\mathfrak l \mid \ell} \bar{\mathcal D}_{\mathfrak l^{v_{\mathfrak l}}}$. Note that, for the same reason as in the proof of Proposition~\ref{prop_local density mod forms bad prime essential multiplicativity}, we may assume $h_{F,f}(\mathfrak o) \neq 0$ and $h_{F,f}(\mathfrak l^{v_{\mathfrak l}} \neq 0$ for any $\mathfrak l \mid \ell$. Under the identification $A_{\mathfrak o}(d) = \prod_{\mathfrak l \mid \ell} A_{\mathfrak l^{v_{\mathfrak l}}}(d_{\mathfrak l})$, $h_{F,f}(\mathfrak o)$ is given by
    \begin{align*}
        h_{F,f}(\mathfrak o)&=\frac{\sum_{d \in \bar{\mathcal D}_{\mathfrak o}} \prod_{\mathfrak l \mid \ell} \#\{g \in A_{\mathfrak l^{v_{\mathfrak l}}}(d_{\mathfrak l}): F(\tr g_{\mathfrak l}) \equiv 0 \pmod{\mathfrak l^{v_{\mathfrak l}}}\}}{\sum_{d \in \bar{\mathcal D}_{\mathfrak o}} \prod_{\mathfrak l \mid \ell} |A_{\mathfrak l^{v_{\mathfrak l}}}(d_{\mathfrak l})|} \\
        &= \frac{\sum_{d \in \bar{\mathcal D}_{\mathfrak o}} \prod_{\mathfrak l \mid \ell} h_F(\mathfrak l^{v_{\mathfrak l}}; \mathcal A_{\mathfrak l}(d_{\mathfrak l})) M_{\mathfrak l^{v_{\mathfrak l}}}}{\sum_{d \in \bar{\mathcal D}_{\mathfrak o}} \prod_{\mathfrak l \mid \ell} M_{\mathfrak l^{v_{\mathfrak l}}}} \\
        &= \frac{1}{|\bar{\mathcal D}_{\mathfrak o}|} \sum_{d \in \bar{\mathcal D}_{\mathfrak o}} \prod_{\mathfrak l \mid \ell} h_F(\mathfrak l^{v_{\mathfrak l}}; \mathcal A_{\mathfrak l}(d_{\mathfrak l})) \\
        &= \frac{1}{|\bar{\mathcal D}_{\mathfrak o}|} \sum_{d \in \bar{\mathcal D}_{\mathfrak o}} \prod_{\mathfrak l \mid \ell} h_{F,f}(\mathfrak l^{v_{\mathfrak l}})(1+\varepsilon_{\mathfrak l^{v_{\mathfrak l}}}(d_{\mathfrak l})) \\
        &= \prod_{\mathfrak l \mid \ell} h_{F,f}(\mathfrak l^{v_{\mathfrak l}}) \bigg(\frac{1}{|\bar{\mathcal D}_{\mathfrak o}|} \sum_{d \in \bar{\mathcal D}_{\mathfrak o}}\prod_{\mathfrak l \mid \ell} (1+\varepsilon_{\mathfrak l^{v_{\mathfrak l}}}(d_{\mathfrak l})) \bigg).
    \end{align*} 
    The properties of $\varepsilon_{\mathfrak l^{v_{\mathfrak l}}}(d_{\mathfrak l})$ established above leads us
    \begin{align*}
        \frac{1}{|\bar{\mathcal D}_{\mathfrak o}|} \sum_{d \in \bar{\mathcal D}_{\mathfrak o}}\prod_{\mathfrak l \mid \ell} (1+\varepsilon_{\mathfrak l^{v_{\mathfrak l}}}(d_{\mathfrak l})) \ll 1 +\frac{1}{|\bar{\mathcal D}_{\mathfrak o}|} \sum_{\substack{J \subset \{ \mathfrak l \mid \ell\} \\ |J| \geq 2}} \prod_{\mathfrak l \in J} q_{\mathfrak l}^{-1} = 1+ O_E(\ell^{-2}).
    \end{align*}
    This completes the proof.
\end{proof}



    The following theorem is a generalization of \cite[Theorem~A.3]{LL26} for $r=1$. We note that \cite[Theorem~A.3]{LL26} is false as stated, although this does not seriously affect the argument for the main result in either \cite{LL26} or \cite{Lee26}; it should be corrected by removing the condition $\ell>q$ from the product.
    
\begin{theorem}\label{thm_h_F,f(mathfrak o) is essentially multiplicative}
    Let
    \begin{align*}
        h_{F,f}(\mathfrak o):=\frac{|\mathcal C_{F,f,\mathfrak o}|}{|\bar{\mathcal A}_{f,\mathfrak o}|}.
    \end{align*}
    Then $h_{F,f}(\mathfrak o)$ is essentially multiplicative in the following sense:
    \begin{align*}
        h_{F,f}(\mathfrak o) \ll \prod_{\mathfrak l^v \mid\mid \mathfrak o}h_{F,f}(\mathfrak l^v)
    \end{align*}
    for all $\mathfrak o$ and
    \begin{align*}
        h_{F,f}(\mathfrak o) \gg \prod_{\mathfrak l^v \mid\mid \mathfrak o}h_{F,f}(\mathfrak l^v)
    \end{align*}
    for $\mathfrak o$ whose prime factor $\mathfrak l$ satisfy $|\kappa(\mathfrak l)|>B$.
\end{theorem}

\begin{proof}
    Let $\mathfrak o_{\ell}:=\prod_{\mathfrak l \mid \ell} \mathfrak l^{v_{\mathfrak l}(\mathfrak o)}$ so that
    \begin{align*}
        \mathfrak o&=\prod_{\ell \mid \Nm \mathfrak o} \mathfrak o_{\ell} = \prod_{\ell \mid \Nm\mathfrak o}\prod_{\mathfrak l \mid \ell}\mathfrak l^{v_{\ell}(\mathfrak o)}.
    \end{align*}
    Let $\rho_{f,\ell}:=\prod_{\mathfrak l \mid \ell} \rho_{f,\mathfrak l}$ and $\bar\rho_{f,\ell}=\prod_{\mathfrak l \mid \ell} \bar\rho_{f,\mathfrak l}$ be its mod $\ell$ reduction.
    According to \cite[Theorems~2.2.2 and 2.3.1]{Loe17}, as $f$ has trivial inner twist, the adelic image is given by
    \begin{align*}
        \rho_f(\Gal_{\Q})=\mathcal U \times \prod_{\ell \not\in S} \mathcal A_{\ell}
    \end{align*}
    for some open subgroup $\mathcal U \subset \mathcal A_S$, where $S$ is a finite set of primes. We may enlarge $S$ if necessary so that any prime that does not belong to $S$ is good. Let $\mathfrak o_{\rm bad}=\prod_{\ell \in S}\mathfrak o_{\ell}$ and $\mathfrak o_{\rm good}=\prod_{\ell \not\in S}\mathfrak o_{\ell}$. Note that $F(\tr g) \equiv 0 \pmod{\mathfrak o}$ if and only if $F(\tr g) \equiv 0 \pmod{\mathfrak o_{\rm bad}}$ and $F(\tr g) \equiv 0 \pmod{\mathfrak o_{\rm good}}$. Thus we have
    \begin{align*}
        h_{F,f}(\mathfrak o)=h_{F,f}(\mathfrak o_{\rm bad})h_{F,f}(\mathfrak o_{\rm good}).
    \end{align*}
    For the bad prime part, by Proposition~\ref{prop_local density mod forms bad prime essential multiplicativity},
    \begin{align*}
        h_{F,f}(\mathfrak o_{\rm bad})=h_F(\mathfrak o_{\rm bad};\mathcal U) \ll_{F,f} \prod_{\mathfrak l \mid \mathfrak o_{\rm bad}}h_F(\mathfrak l^{v_{\mathfrak l}(\mathfrak o)};\mathcal U)=\prod_{\mathfrak l^v \mid\mid \mathfrak o_{\rm bad}} h_{F,f}(\mathfrak l^v).
    \end{align*}
    Clearly, one also has
    \begin{align*}
        h_{F,f}(\mathfrak o_{\rm good})=h_{F,f}\bigg(\prod_{\substack{\ell \mid \Nm \mathfrak o \\ \ell \gg_{F,f} 1}} \mathfrak o_{\ell} \bigg)  \ll \prod_{\substack{\ell \mid \Nm \mathfrak o \\ \ell \gg_{F,f} 1}}h_{F,f}(\mathfrak o_{\ell})
    \end{align*}
    for all $\mathfrak o$ and
    \begin{align*}
        h_{F,f}(\mathfrak o) \gg \prod_{\ell \mid \Nm \mathfrak o}h_{F,f}(\mathfrak o_{\ell})
    \end{align*}
    for $\mathfrak o$ whose prime factors $\mathfrak l$ satisfy $q_{\mathfrak l} \gg_{F,f} 1$.

    By Proposition~\ref{prop_local density mod forms one good prime decomposition, stronger essetinal multiplicativity}, for $\ell \gg 1$,
    \begin{align*}
        h_{F,f}(\mathfrak o_{\ell}) = \bigg(1+O\bigg(\frac{1}{\ell^2}\bigg)\bigg)\prod_{\substack{\mathfrak l^v \mid\mid \mathfrak o_{\ell}}} h_{F,f}(\mathfrak l^v).
    \end{align*}
    Therefore, we obtain
    \begin{align*}
        h_{F,f}(\mathfrak o_{\rm good})\ll \prod_{\substack{\ell \mid \Nm \mathfrak o \\ \ell \gg 1}} h_{F,f}(\mathfrak o_{\ell})&=\prod_{\substack{\ell \mid \Nm \mathfrak o \\ \ell \gg 1}}\bigg(1+O\bigg(\frac{1}{\ell^2}\bigg)\bigg)\prod_{\substack{\mathfrak l^v \mid\mid \mathfrak o_{\ell}}} h_{F,f}(\mathfrak l^v) \\
        & \ll \prod_{\ell}\bigg(1+O\bigg(\frac{1}{\ell^2}\bigg)\bigg)\times \prod_{\substack{\mathfrak l^v \mid\mid \mathfrak o \\ \ell \gg_{F,f} 1}} h_{F,f}(\mathfrak l^v) \\
        &\ll \prod_{\substack{\mathfrak l^v \mid\mid \mathfrak o \\ \ell \gg_{F,f} 1}} h_{F,f}(\mathfrak l^v).
    \end{align*}
    The second assertion follows in a similar way.
\end{proof}

\begin{lemma}\label{lem_C_f,F,lambda^v is nonempty}\footnote{This lemma was suggested by ChatGPT to resolve the gap concerning bad primes in the proof of Theorem~\ref{thm_application 2: modular forms} below.}
    Let $\mathfrak l$ be a finite place of $E$ lying over a rational prime $\ell$. If $\mathfrak l^v \mid F(a_f(\ell))$ for $v \geq 1$ and $\mathfrak l \gg_f 1$, then $\mathcal C_{f,F,\mathfrak l^v}$ is nonempty.
\end{lemma}

\begin{proof}
    From $\mathfrak l \gg_f 1$, we have (cf. \cite[Theorem~2.2.2]{Loe17})
    \begin{align*}
        \rho_{f,\mathfrak l}(\Gal_{\Q})=\{g \in \GL_2(\mathcal O_{E,\mathfrak l}): \det g \in (\Z_{\ell}^{\times})^{k-1}\}.
    \end{align*}
    Thus, for any $d \in (\Z_{\ell}^{\times})^{k-1} \subset (O_{E,\mathfrak l}^{\times})^{k-1}$ whose reduction modulo $\mathfrak l^v$ is also denoted by $d$, the matrix
    \begin{align*}
        g:=\begin{pmatrix}
            a_f(\ell) & d \\ -1 & 0
        \end{pmatrix} \in \mathcal A_{f,\mathfrak l^v}
    \end{align*}
    satisfies $F(a_f(\ell))=F(\tr(g)) \equiv 0 \pmod{\mathfrak l^v}$. This means $g \in \mathcal C_{f,F,\mathfrak l^v}$.
\end{proof}

Now we prove Theorem~\ref{thm_application 2: modular forms}.

\begin{proof}[Proof of Theorem~\ref{thm_application 2: modular forms}]
    Denote by $K_{\mathfrak o}=\Q(\bar\rho_{f,\mathfrak o})=\overline \Q^{\ker \bar\rho_{f,\mathfrak o}}$ the field fixed by the kernel of $\bar\rho_{f, \mathfrak o}$. We have $\Gal(K_{\mathfrak o}/\Q) \simeq A_{f,\mathfrak o}$. The discriminant of $K_{\mathfrak o} / \Q$ will be denoted by $\Delta_{\mathfrak o}$. Let
    \begin{align*}
        \Pi^*(x,I;f):=\{p\leq x : p\nmid N, \ \theta_f(p) \in I\}.
    \end{align*}
    By the effective Chebotarev theorem (\cite[Theorem~1.1]{LO77}) under GRH,
    \begin{align*}
        \#\{p\in \Pi^*(x,[0,\pi];f): \mathfrak o \mid F(a_f(p))\}&=\frac{|\mathcal C_{F,f,\mathfrak o}|}{|\mathcal A_{f,\mathfrak o}|}\pi(x)+O\bigg(\frac{|\mathcal C_{F,f,\mathfrak o}|}{|\mathcal A_{f,\mathfrak o}|}x^{1/2}\log(|\Delta_{\mathfrak o}|x^{[K_{\mathfrak o}:\Q]})+\log|\Delta_{\mathfrak o}|\bigg) \\
        &=h_{F,f}(\mathfrak o)\pi(x)+O\bigg(h_{F,f}(\mathfrak o)x^{1/2}\log(|\Delta_{\mathfrak o}|x^{[K_{\mathfrak o}:\Q]})+\log|\Delta_{\mathfrak o}|\bigg)
    \end{align*}
    We apply the same argument as in \cite[p.72--73]{MM84} as follows. We have
    \begin{align*}
        [K_{\mathfrak o}:\Q]=|\mathcal A_{f,\mathfrak o}| \leq |\GL_2(\mathcal O_E/\mathfrak o)| \leq (\Nm \mathfrak o)^4
    \end{align*}
    and by \cite[\S~1.4]{Ser81},
    \begin{align*}
        \log|\Delta_{\mathfrak o}| &\leq [K_{\mathfrak o}:\Q]\log((\Nm \mathfrak o)N[K_{\mathfrak o}:\Q]) \\
        & \leq (\Nm \mathfrak o)^4 \log((\Nm \mathfrak o)^5 N).
    \end{align*}
    Consequently,
    \begin{align*}
        \#\{p\in \Pi^*(x,[0,\pi];f): \mathfrak o \mid F(a_f(p))\}=h_{F,f}(\mathfrak o)\pi(x)+O(x^{1/2}(\Nm \mathfrak o)^4\log(xN\Nm \mathfrak o)).
    \end{align*}
    According to \cite[\S~7.2]{Ser81}, under GRH, we have
    \begin{align}\label{eqn_F(a_f(p))=0 exclude}
        \#\{p\leq x: F(a_f(p))= 0\} \ll_{F,f} x^{7/8}(\log x)^{1/2}.
    \end{align}
    Thus we have
    \begin{align*}
        \#\{p\in \Pi(x,[0,\pi];F,f): \mathfrak o \mid F(a_f(p))\}=h_{F,f}(\mathfrak o)\pi(x)+O_{F,f}((\Nm \mathfrak o)^4\log \Nm \mathfrak o +x^{7/8}\log x).
    \end{align*}

    When $I\neq [0,\pi]$, under GRH and \cite[Hypothesis~A]{LL26}, we can also obtain a similar estimate. Following the same argument as in \cite[\S~6.2]{LL26}, we use the effective Chebotarev--Sato--Tate theorem (\cite[Theorem~1.3]{LL26}, see \cite[Theorem~1.2]{Won19} for a general statement on Hilbert modular forms) to obtain
    \begin{align*}
        &\#\{p\in \Pi^*(x,I;f): \mathfrak o \mid F(a_f(p))\} \\ &=h_{F,f}(\mathfrak o)\mu_{ST}(I)\pi(x) +O\bigl(x^{3/4}|\mathcal C_{F,f,\mathfrak o}|(\log(Nx[K_{\mathfrak o}:\Q]|\Delta_{\mathfrak o}|\prod_{p \ \text{ramified in} \ K_{\mathfrak o}}p))^{1/2}\bigl) \\ 
        &=h_{F,f}(\mathfrak o)\mu_{ST}(I)\pi(x) + O_f\bigl(x^{3/4}|\mathcal C_{F,f,\mathfrak o}|(\log x)^{1/2}+x^{3/4}|\mathcal C_{F,f,\mathfrak o}|(\Nm \mathfrak o)^2 (\log \Nm \mathfrak o)^{1/2}\bigl) \\
        &=h_{F,f}(\mathfrak o)\mu_{ST}(I)\pi(x) + O_f\bigl(x^{3/4}(\log x)^{1/2}(\Nm \mathfrak o)^{6+\epsilon}\bigl).
    \end{align*}
    By~\eqref{eqn_F(a_f(p))=0 exclude},
    \begin{align*}
        \#\{p\in \Pi(x,I;F,f): \mathfrak o \mid F(a_f(p))\} = h_{F,f}(\mathfrak o)\mu_{ST}(I)\pi(x) + O_{F,f}\bigl(x^{7/8}(\log x)^{1/2}(\Nm \mathfrak o)^{6+\epsilon}\bigl).
    \end{align*}
    Under GRH for the symmetric power $L$-functions of modular forms, we have (\cite[Theorem~1.3]{Tho21}) 
    \begin{align*}
        |\Pi(x,I;F,f)_{\rm good}| - \mu_{ST}(I)\pi(x) \ll_{F,f} x^{7/8} (\log x)^{1/2},
    \end{align*}
    and therefore,
    \begin{align*}
        h(x,\mathfrak o)=h_{F,f}(\mathfrak o) + O_{F,f}\bigl(\mu_{ST}(I)^{-1}x^{-1/8}(\log x)^{3/2}(\Nm \mathfrak o)^{6+\epsilon}\bigl).
    \end{align*}  
    This is the modular analogue of~\eqref{eqn_h(x,delta) and N_F(delta) polynomial}. We prove that $h(x,D)$ is essentially multiplicative in $D$ with $\Nm D \leq x^c$ for a suitable $c \in (0,1)$. Note that the condition $\mathcal C_{F,f,\mathfrak o}=\emptyset$ implies $h(x,\mathfrak o)=0$. Indeed, if $\mathcal C_{f,F,\mathfrak l^v}=\emptyset$, then no $p \gg_f 1$ with $p \neq \ell$ satisfies 
    \begin{align*}
        \mathfrak l^v \mid F(a_f(p))=F(\tr(\Frob_p)). 
    \end{align*}
    If $p=\ell$, then Lemma~\ref{lem_C_f,F,lambda^v is nonempty} gives the same conclusion. Therefore, no $p \gg_f 1$ satisfies $\mathfrak l^v \mid F(a_f(p))$, and hence no $p \gg_f 1$ satisfies $\mathfrak o \mid F(a_f(p))$.

    We may assume $\mathcal C_{F,f,\mathfrak o}\neq \emptyset$. Since
    \begin{align*}
        h_{F,f}(\mathfrak o) \geq \frac{1}{|\mathcal A_{f,\mathfrak o}|} \geq \frac{1}{(\Nm \mathfrak o)^4},
    \end{align*}
    one has
    \begin{align*}
        \frac{h(x,\mathfrak o)}{h_{F,f}(\mathfrak o)}&=1+O_{F,f}\bigg(\frac{\mu_{ST}(I)^{-1}x^{-1/8}(\log x)^{1/2}(\Nm \mathfrak o)^{6+\epsilon}}{h_{F,f}(\mathfrak o)}\bigg) \\
        &=1+O_{F,f}(\mu_{ST}(I)^{-1}x^{-1/8}(\log x)^{1/2}(\Nm \mathfrak o)^{10+\epsilon}).
    \end{align*}
    The big $O$-term is $o(1)$ if we take $c>0$ sufficiently small, specifically so that
    \begin{align*}
        c<\frac{1}{8(10+\epsilon)}.
    \end{align*}
    Essential multiplicativity of $h(x,D)$ follows from Theorem~\ref{thm_h_F,f(mathfrak o) is essentially multiplicative}.

    If $y^v \leq x^c$, for $\ell \gg 1$, we have
    \begin{align*}
        h_{F,f}(\mathfrak l^v)=\sum_{i=1}^s h_{F_i,f}(\mathfrak l^{\lceil v/e_i \rceil}).
    \end{align*}
    If $1\leq v \leq e_i$, then $h_{F_i,f}(\mathfrak l^{\lceil v/e_i \rceil})= h_{F_i,f}(\mathfrak l)=N_{F_i}(\mathfrak l)q_{\mathfrak l}^{-1}(1+O(q^{-1}))$ for $\ell \gg 1$ by Lemma~\ref{lem_A lambda^n d, t number}, where
    \begin{align*}
        N_{F_i}(\mathfrak l):=\#\{t \in \kappa(\mathfrak l): F_i(t)=0\}
    \end{align*}
    as in \S\S~\ref{subsec_polynomial_values}. If $v>e_i$, then $h_{F_i,f}(\mathfrak l^{\lceil v/e_i \rceil}) \ll q_{\mathfrak l}^{-2}$. Thus we have
    \begin{align*}
        \sum_{q_{\mathfrak l} \leq y}h_{F_i,f}(\mathfrak l^v)= \sum_{q_{\mathfrak l} \leq y}\frac{N_{F_i}(\mathfrak l)}{q_{\mathfrak l}}+O_F(1)
    \end{align*}
    and hence
    \begin{align*}
        \sum_{q_{\mathfrak l} \leq y}h_{F,f}(\mathfrak l^v)&=\sum_{i : v\leq e_i}\sum_{q_{\mathfrak l} \leq y}\frac{N_{F_i}(\mathfrak l)}{q_{\mathfrak l}}+O_F(1) \\
        &= \sum_{i : v\leq e_i} \log\log y + O_F(1) \\
        &= \kappa_v \log\log y + O_F(1).
    \end{align*}
    On the other hand, as we have shown above,
    \begin{align*}
        h(x,\mathfrak l^n)-h_{F,f}(\mathfrak l^n) \ll_{F,f} \mu_{ST}(I)^{-1}x^{-1/8}(\log x)^{3/2}q_{\mathfrak l}^{n(6+\epsilon)}
    \end{align*}
    and hence
    \begin{align*}
        \sum_{q_{\mathfrak l} \leq y} |h(x,\mathfrak l^n)-h_{F,f}(\mathfrak l^n)| = o(1).
    \end{align*}
    The remaining steps are the same as in the proof of Theorem~\ref{thm_application 1: polynomial value}.
\end{proof}

\subsection{Polynomial values over function fields}

In this subsection, we prove Theorem~\ref{thm_application 3: function field polynomial}. Henceforth, we adopt the notation in Theorems~\ref{thm_application 3: function field polynomial_weak} and \ref{thm_application 3: function field polynomial}. We set $X=\Spec R$.

Let $N_F(\delta)$ be the number of roots of $F$ modulo $\delta$, i.e., the number of $u \in R/\delta R$ such that $F(u)=0 \in R/\delta R$. By the Chinese Remainder Theorem, it is multiplicative in $\delta$. 

\begin{lemma}\label{lem_inseparability_zero}
    If $F \in R[T]$ is irreducible and inseparable, then $N_F(\ell^e)=0$ for $e \geq 2$.
\end{lemma}

\begin{proof}
    Suppose that there exists $b \in R$ with $\ell^e \mid F(b)$ for $e\geq 2$. Since $F(b)$ can be seen as a polynomial in $U$, we have $\ell \mid (\partial_U F)(b)$; in particular $b$ is a common root of $F(T)$ and $(\partial_U F)(T)$ mod $\ell$.

    On the other hand, since $F$ is irreducible in $T$, either $\gcd_T(F, \partial_U F)=1$ or $\partial_U F=0$. As we have a common root $b$, this implies $\partial_U F=0$. As $F$ is inseparable, we also have $\partial_T F=0$, so
    \begin{align*}
        F \in \F_q[U^{p},T^{p}].
    \end{align*}
    In particular, this implies that $F$ is the $p$th power of some polynomial in $\F_q[U,T]=R[T]$, which yields a contradiction.
\end{proof}

\begin{lemma}\label{lem_sieve_density_sequence_function_field_case}
    We have
    \begin{align*}
        \sum_{\Nm \ell^v \leq y} \frac{N_F(\ell^v)}{\Nm \ell^v} = \kappa_v \log\log y + O_{q,F}(1),
    \end{align*}
    where $\ell$ ranges over the primes in $R$.
\end{lemma}

\begin{proof}
    First, we prove the assertion when $F$ is irreducible and $v=1$. Let $K_F=K[T]/(F)$ be a finite extension of $K$, and let $\F_{q^{c_F}}$ be its constant field. Then the norm $\Nm \mathfrak l$ of a prime element $\mathfrak l$ of $K_F$ is given by $q^{c_F \deg \mathfrak l}$. Denote by $\mathcal O_F$ the ring of integers of $K_F$.

    Note that, according to the Dedekind--Kummer theorem, $N_F(\ell)$ equals the number of primes $\mathfrak l$ of $K_F$ lying over $\ell$ for which the inertia degree $
        \mathfrak f(\mathfrak l \mid \ell):=[\mathcal O_F/\mathfrak l: R/\ell]
    $
    is $1$. Since $\Nm \mathfrak l = (\Nm \ell)^{\mathfrak f(\mathfrak l \mid \ell)}$, we have
    \begin{align*}
        \sum_{\Nm \ell\leq y}\frac{N_{F}(\ell)}{\Nm \ell}= \sum_{\Nm \ell \leq y}\sum_{\substack{\mathfrak l \in \cl(\Spec \mathcal O_F) \\ \mathfrak l \mid \ell \\ \mathfrak f(\mathfrak l \mid \ell)=1}} \frac{1}{\Nm \mathfrak l}= \sum_{\substack{\mathfrak l \in \cl(\Spec \mathcal O_F) \\ \Nm \mathfrak l \leq y}} \frac{1}{\Nm \mathfrak l}+O_{q,F}(1).
    \end{align*}
    The last equality follows from the fact that the sum of $1/\Nm\mathfrak l=1/(\Nm \ell)^{\mathfrak f(\mathfrak l \mid \ell)}$ for all $\mathfrak l$ with $\mathfrak f(\mathfrak l \mid \ell) \geq 2$ converges to a constant depending on $q$ and $F$. 
    
    By Rosen~\cite[Theorem~5.12]{Ros02},
    \begin{align*}
        \sum_{\substack{\mathfrak l \in \cl(\Spec \mathcal O_F) \\ \Nm \mathfrak l \leq y}} \frac{1}{\Nm \mathfrak l} &= \sum_{r=1}^{\lfloor \log_{q^{c_F}}\rfloor} \sum_{\substack{\mathfrak l \in \cl(\Spec \mathcal O_F) \\\deg \mathfrak l=r}} \frac{1}{\Nm \mathfrak l} = \sum_{r=1}^{\lfloor \log_{q^{c_F}}\rfloor} \#\{\mathfrak l \in \cl(\Spec \mathcal O_F): \deg \mathfrak l=r\} \frac{1}{q^{c_F r}} \\
        &=\sum_{r=1}^{\lfloor \log_{q^{c_F}}\rfloor} \frac{1}{q^{c_F r}}\bigg(\frac{q^{c_F r}}{r}+O_{K_i}\bigg(\frac{q^{c_F r/2}}{r}\bigg)\bigg) = \log\log y + O_{q,F}(1).
    \end{align*}
    We now consider the general case $F=c\prod_{i=1}^s F_i^{e_i}$ with $v \geq 1$. For all but finitely many $\ell$, $F(b)$ is divisible by $\ell^v$ if and only if $F_i(b)$ is divisible by $\ell^{\lceil v/e_i\rceil}$ for some $i$. If $e_i \geq v$, then $N_{F_i}(\ell^{\lceil v/e_i \rceil})/\Nm \ell^{\lceil v/e_i \rceil} = N_{F_i}(\ell)/\Nm \ell$. If $e_i < v$, then by Lemma~\ref{lem_inseparability_zero} we may assume $F_i$ is separable, and hence by Hensel's lemma, $N_{F_i}(\ell^{\lceil v/e_i \rceil})/\Nm \ell^{\lceil v/e_i \rceil} \ll_{F_i} (\Nm \ell)^{-2}$. Its contribution to the sum is $O_{q,F}(1)$. Consequently, 
    \begin{align*}
        \sum_{\Nm \ell^v \leq y} \frac{N_F(\ell^v)}{\Nm \ell^v} = \sum_{i: e_i \geq v} \sum_{\Nm \ell \leq y} \frac{N_{F_i}(\ell)}{\Nm \ell} + O_{q,F}(1).
    \end{align*}
    We complete the proof by applying the assertion to each irreducible polynomial $F_i$.
\end{proof}

\begin{proof}[Proof of Theorem~\ref{thm_application 3: function field polynomial}]
    Let $n_0$ be a positive integer depending on $F$ such that if $r \geq n_0$, then there is no $b \in S_r$ with $F(b)=0$ so that $|S_r|=q^r$. For $x \in (q^{n_0}, \infty)$, set
    \begin{align*}
        A_x=\bigsqcup_{r=n_0}^{\lfloor \log_q x \rfloor} (S_r \times \{1,2,\ldots,q^{r^3}\}).
    \end{align*}
    Clearly, $A_x \subseteq A_y$ if $x\leq y$. We also define $\n(b,j)$ to be the principal divisor $D_{F(b)}$ associated with $F(b) \in R$ for $(b,j) \in A$.
    Then by definition one has
    \begin{align*}
        \pi(x,D_{\delta})=\#\{(b,j) \in A_x: \delta \mid F(b)\} = \sum_{r=n_0}^{\lfloor \log_q x \rfloor} q^{r^3}|S_r|h_{F,r}(\delta)
    \end{align*}
    and consequently,
    \begin{align*}
        h(x,D_{\delta})=\frac{\sum_{r=n_0}^{\lfloor \log_q x \rfloor} q^{r^3}|S_r|h_{F,r}(\delta)}{\sum_{r=n_0}^{\lfloor \log_q x \rfloor} q^{r^3}|S_r|}.
    \end{align*}
    Note that if $n \geq \deg \delta$ and $n\geq n_0$, then the number of $b \in S_n$ with $\delta \mid F(b)$ is $N_F(M)q^{n-\deg \delta}$. Thus
    \begin{align*}
        h_{F,n}(\delta)=\frac{N_F(\delta)q^{n-\deg \delta}}{q^n}=\frac{N_F(\delta)}{\Nm \delta},
    \end{align*}
    which is independent of the choice of $n$.

    As $h_{F,r}(\delta) \in [0,1]$ for all $r \in \N$ and $\delta \in R$, we have
    \begin{align*}
        |h(q^n,D_{\delta})-h_{F,n}(\delta)|&=\left|\frac{\sum_{r=n_0}^{\deg \delta-1} q^{r^3}|S_r|(h_{F,r}(\delta)-h_{F,n}(\delta))}{\sum_{r=n_0}^{n} q^{r^3}|S_r|}\right| \leq \frac{\sum_{r=n_0}^{n-1} q^{r^3+r}}{\sum_{r=n_0}^{n} q^{r^3+r}} \\
        &\leq \frac{\sum_{r=n_0}^{n-1} q^{r^3+r}}{q^{n^3+n}} \leq \frac{nq^{(n-1)^3+(n-1)}}{q^{n^3+n}} \ll_q q^{-n^2}.
    \end{align*}
    Fix any $c \in (0,1)$. If $n=\lfloor \log_q x \rfloor$ and $\Nm \delta \leq x^c$, then
    \begin{align*}
        h(x,D_{\delta})&=h_{F,n}(\delta)(1+O_q((\Nm \delta)q^{-n^2}) \\
        &=\frac{N_F(\delta)}{\Nm \delta}(1+O_q(x^cq^{-n^2})).
    \end{align*}
    By the multiplicativity of $\frac{N_F(\delta)}{\Nm \delta}$, one has
    \begin{align*}
        \prod_{\ell^v \mid\mid \delta} h(x,D_{\ell^v})&=\prod_{\ell^v \mid\mid \delta}\frac{N_F(\ell^v)}{\Nm \ell^v}(1+O_q(x^cq^{-n^2}))  \\
        &=\frac{N_F(\delta)}{\Nm \delta}(1+O_q(x^c(\log x)q^{-n^2})) \\
        &= h(x,D_{\delta})(1+o_q(1)),
    \end{align*}
    which implies the essential multiplicativity of $h(x,D_{\delta})$ in $D_{\delta}$ with $\Nm \delta \leq x^c$.

    According to Lemma~\ref{lem_sieve_density_sequence_function_field_case}, we have
    \begin{align*}
        \sum_{\Nm \ell^v \leq x} h(x,v[\ell]) &= \sum_{\Nm \ell^v \leq x} \frac{N_F(\ell^v)}{\Nm \ell^v} (1+O_q(x^c q^{-n^2})) \\
        &= (\kappa_v \log\log x + O_{q,F}(1))(1+O_q(x^c q^{-n^2})).
    \end{align*}
    As $x^c q^{-n^2}\log\log x \to 0$ as $x \to \infty$, $h$ is of order~$1$ with the sieve dimension sequence $(\kappa_v)_{v \geq 1}$. The dimension sequence clearly has finite length. If we consider
    \begin{align*}
        \sum_x \sum_{\ell \in \cl(X)} h(x,v[\ell])^2,
    \end{align*}
    then, by an argument similar to the one above, it is
    \begin{align*}
        =\sup_x \sum_{\ell \in \cl(X)} \frac{N_F(\ell^v)^2}{\Nm \ell^{2v}}(1+O_q(x^c q^{-n^2}))\ll_F \sup_x \sum_{\ell \in \cl(X)} \frac{1}{\Nm \ell^{2v}} < \infty.
    \end{align*}
    Lastly, for all but finitely many $\ell$, we have
    \begin{align*}
        h(x,v[\ell]) \ll_F \sum_{i=1}^s |\kappa(\ell)|^{-\lceil v/e_i \rceil} \ll_F |\kappa(\ell)|^{-v/E}
    \end{align*}
    by an argument similar to the one above, so \eqref{eqn_h(x,vp) for v>E bound} holds for $B=0$ and $\alpha=1/E$. Therefore, the conditions (i)--(iii) of Theorem~\ref{thm_main_geometric} are satisfied, so the result follows.
\end{proof}

\end{document}